\documentclass[11pt]{article}
\usepackage{times,color}
\usepackage[top=2cm,bottom=2cm,left=2cm,right=2cm]{geometry}
\usepackage{amssymb, graphicx, amsmath,bm}
\usepackage{amsthm}
\usepackage{subfigure}
\usepackage{algorithm}
\usepackage{algpseudocode}
\usepackage{cases}
 \usepackage{epstopdf}
\usepackage{epsfig}

\usepackage[T1]{fontenc}
\usepackage[utf8]{inputenc}
\usepackage[colorlinks,citecolor=red,linkcolor=blue]{hyperref}
\newtheorem{theorem}{Theorem}[]

\newtheorem{remark}[theorem]{Remark}

\newtheorem{example}{Example}
\makeatletter
\@addtoreset{equation}{section}
\makeatother

\floatname{algorithm}{Subroutine}
\usepackage{tikz}
\usetikzlibrary{shapes.geometric, arrows}
\tikzstyle{start} = [rectangle, rounded corners, draw=black, fill=red!30]
\tikzstyle{rectangle_rounded_corners} = [rectangle, rounded corners, text centered, draw=black, fill=red!30]
\tikzstyle{rectangle_rounded_green} = [rectangle, rounded corners, text width=8cm, text centered, draw=black, fill=green!30]
\tikzstyle{rectangle_out} = [rectangle, rounded corners, text centered, draw=black, fill=green!30]
\tikzstyle{trapezium_blue} = [trapezium, trapezium left angle=70, trapezium right angle=110, text centered, draw=black, fill=blue!30]
\tikzstyle{trapezium_orange} = [trapezium, trapezium left angle=70, trapezium right angle=110, text centered, draw=black, fill=orange!30]
\tikzstyle{rec} = [rectangle,  text centered, draw=black, fill=orange!30]
\tikzstyle{rec_lable} = [rectangle,  text centered]
\tikzstyle{rec_long} = [rectangle, text width=7cm, text centered, draw=black, fill=orange!30]
\tikzstyle{ellipse_red} = [ellipse, text centered, draw=black, fill=red!30]
\tikzstyle{arrow} = [thick,->,>=stealth]

\title{An ADMM-Newton-CNN Numerical Approach to a TV Model for Identifying Discontinuous Diffusion Coefficients in Elliptic Equations: Convex Case with Gradient Observations}
\author{Wenyi Tian\thanks{Center for Applied Mathematics, Tianjin University, Tianjin 300072, China. This author was partially supported by the National Natural Science Foundation of China (No. 12071343). Email: twymath@gmail.com}
        \and
        Xiaoming Yuan\thanks{Department of Mathematics, The University of Hong Kong, Hong Kong, China. Email: xmyuan@hku.hk}
        \and
        Hangrui Yue\thanks{Corresponding author. Department of Mathematics, The University of Hong Kong, Hong Kong, China. Email: yuehangrui@gmail.com}
        }
\begin{document}
\maketitle
\begin{abstract}
Identifying the discontinuous diffusion coefficient in an elliptic equation with observation data of the gradient of the solution is an important nonlinear and ill-posed inverse problem. Models with total variational (TV) regularization have been widely studied for this problem, while the theoretically required nonsmoothness property of the TV regularization and the hidden convexity of the models are usually sacrificed when numerical schemes are considered in the literature. In this paper, we show that the favorable nonsmoothness and convexity properties can be entirely kept if the well-known alternating direction method of multipliers (ADMM) is applied to the TV-regularized models, hence it is meaningful to consider designing numerical schemes based on the ADMM. Moreover, we show that one of the ADMM subproblems can be well solved by the active-set Newton method along with the Schur complement reduction method, and the other one can be efficiently solved by the deep convolutional neural network (CNN). The resulting ADMM-Newton-CNN approach is demonstrated to be easily implementable and very efficient even for higher-dimensional spaces with fine mesh discretization.
\end{abstract}

\section{Introduction}
Consider the canonical elliptic equation
\begin{equation}\label{eq:IDCEE}
  \left\{
  \begin{aligned}
    &-\nabla\cdot(q(x)\nabla u(x))=f(x), &&x\in\Omega,\\
    &u(x)=0, &&x\in\Gamma,
  \end{aligned}\right.
\end{equation}
where $\Omega$ is a bounded polyhedral domain in $\mathbb{R}^d$ ($d=1,2,3$) with a piecewise smooth boundary $\Gamma:=\partial\Omega$; $u:\Omega\cup\Gamma\rightarrow\mathbb{R}$ with $u\in H_0^1(\Omega)$; $q:\Omega\rightarrow\mathbb{R}$ with $q\in L^{\infty}(\Omega)$; and $f:\Omega\rightarrow\mathbb{R}$ with $f\in H^{-1}(\Omega)$ is given. The elliptic equation \eqref{eq:IDCEE} describes various physical phenomena such as the flow of a fluid through some medium with the permeability $q(x)$ and the concentration $u(x)$, and the heat transfer in a material with the conductivity $q(x)$ and the temperature $u(x)$. For the diffusion coefficient $q(x)$, it is often impractical to measure it directly (e.g., when it is the conductivity of a medium), but it is easier to observe the solution $u$ of \eqref{eq:IDCEE} or its gradient \cite{ChenZ:1999,GuentherHMG:1985}. Hence, it is interesting to consider the inverse problem of identifying the diffusion coefficient $q(x)$ with observation data of the solution $u$ of the elliptic equation \eqref{eq:IDCEE} or its gradient. This inverse problem finds applications in various industrial areas such as reservoir simulations, underground water investigations, geophysics and electrical impedance tomography. We refer to the monographs \cite{BanksK:1989,EnglHN:1996} for more introductions. Note that the elliptic equation \eqref{eq:IDCEE} is linear if $q(x)$ is known, but the inverse problem of identifying $q(x)$ is nonlinear.
Also, as mentioned in \cite{EnglHN:1996,KohnL:1988}, $q(x)$ cannot be uniquely determined by $u(x)$ since $q(x)$ can be arbitrary when $u$ is constant on some open subset of $\Omega$. Thus, identifying the diffusion coefficient $q(x)$ of (\ref{eq:IDCEE}) is an ill-posed inverse problem.

\subsection{TV model}
Let us consider the case where observation data of the gradient of the solution $u$ of (\ref{eq:IDCEE}) are available subject to some noise with the noisy level $\delta>0$; it is denoted by $\nabla u_{\delta}\in (L^2(\Omega))^d$. In some literatures such as \cite{ChenZ:1999,HaoQ:2011,Zou:1998}, it has been proposed to recover the discontinuous coefficient $q(x)$ in \eqref{eq:IDCEE} with $\nabla u_{\delta}\in (L^2(\Omega))^d$ via the model
\begin{equation}\label{eq:ROTV}
  \begin{aligned}
    \min_{q,u}&~\Big\{\frac{1}{2}\int_\Omega q|\nabla u-\nabla u_{\delta}|^2\mathrm{d}x+\alpha\int_{\Omega}|\nabla q|\Big\},\\
    \text{s.t.}~&~-\nabla\cdot(q\nabla u)=f,~~(q,u)\in K\times H_0^1(\Omega),
  \end{aligned}
\end{equation}
in which $ \frac{1}{2}\int_\Omega q|\nabla u-\nabla u_{\delta}|^2\mathrm{d}x$ is a data-fidelity term and $\int_{\Omega}|\nabla q|$ is the total variation (TV) regularization term defined in \cite{AttouchBM:2006}. That is, we have
\begin{equation}\label{eq:TV}
  \int_{\Omega}|\nabla q|:=\sup\Big\{\int_{\Omega}q~\mathrm{div} \varphi~\mathrm{d}x: \varphi\in C_c^{1}(\Omega;\mathbb{R}^d),\|\varphi\|_{\infty}\le1\Big\}
\end{equation}
with $\|\varphi\|_{\infty}=\sup_{x\in\Omega}(\sum_{i=1}^d|\varphi_i(x)|^2)^{1/2}$, ``div'' denotes the divergence operator, and $C_c^{1}(\Omega;\mathbb{R}^d)$ is the set of once
continuously differentiable $\mathbb{R}^d$-valued functions with compact support in $\Omega$, see, e.g., \cite{AttouchBM:2006,Ziemer:1989} for more details. Moreover, the admissible set $K$ is
\begin{equation}\label{eq:K}
  K:=\{q\in L^{\infty}\cap BV(\Omega):0<a_0\le q(x)\le a_1, \mbox{~a.e.~in~}\Omega\},
\end{equation}
and the $BV(\Omega)$ space endowed with the norm $\|q\|_{BV}:=\|q\|_{L^1(\Omega)}+\int_{\Omega}|\nabla q|$ is a Banach space; see, e.g., \cite{AttouchBM:2006,Ziemer:1989} for more details. Also, in (\ref{eq:ROTV}), $\alpha>0$ is a parameter determining the relative weights of the data-fidelity and TV regularization terms in the objective functional.

Note that the elliptic equation \eqref{eq:IDCEE} is the Euler-Lagrange equation of the energy functional
$$\displaystyle\frac{1}{2}\int_\Omega(q|\nabla u|^2-2fu)\mathrm{d}x,
$$
and the following identity holds (see, e.g., \cite{Knowles:2001}):
\begin{equation*}
  \frac{1}{2}\int_\Omega q|\nabla u-\nabla u_{\delta}|^2\mathrm{d}x
  =\frac{1}{2}\int_\Omega(q|\nabla u_{\delta}|^2-2fu_{\delta})\mathrm{d}x-\frac{1}{2}\int_\Omega(q|\nabla u|^2-2fu)\mathrm{d}x,
\end{equation*}
where $u_{\delta}$ denotes the approximation to the solution $u$ subject to the noise level $\delta>0$. Hence, the data-fidelity term in \eqref{eq:ROTV} measures the difference of the energy functional of the elliptic equation \eqref{eq:IDCEE} at $u_{\delta}$ and $u$, and it has been widely used in the literature, see, e.g., \cite{ChenZ:1999,HaoQ:2010,HaoQ:2011,HinzeQ:2016,Knowles:2001,KohnV:1984,KohnV:1987,KohnL:1988}. For the TV regularization term, it is capable of reserving the piecewise-constant property and it has found various applications such as image denoising or reconstruction. The TV regularization has also been considered for identifying the diffusion coefficient $q(x)$ of \eqref{eq:IDCEE} because it is generally discontinuous and also owns the piecewise-constant property for many applications such as reservoir simulations and electrical impedance tomography. We refer to, e.g., \cite{ChanT:1997,ChanT:2003,ChenZ:1999,HaoQ:2011}, for more discussions.
One interesting fact is that, as proved in \cite{HaoQ:2010}, although the data-fidelity functional in \eqref{eq:ROTV} is nonconvex with respect to $q$ and $u$ jointly, it is convex with respect to $q$ if $u$ is represented as a function of $q$.

\subsection{ALM for the smoothing $TV_{\epsilon}$ model}

The TV term $\int_{\Omega}|\nabla q|$ in \eqref{eq:ROTV} is not differentiable and it could be difficult to tackle the nonsmoothness property for algorithmic design. In earlier literatures, it is popular to consider smoothing the TV term and then use the smoothing surrogate to replace the original TV term. For instance, in \cite{AcarV:1994,ChanT:1997,ChanT:2003,ChenZ:1999,KeungZ:1998,TaiFEMC:1998}, it is suggested to replace the TV term $\int_{\Omega}|\nabla q|$ with the surrogate
$$
  TV_{\epsilon}(q):=\int_{\Omega}\sqrt{|\nabla q|^2+\epsilon}~\mathrm{d}x,
$$
where $\epsilon>0$ is a smoothing parameter such that $TV_{\epsilon}(q)\rightarrow\int_{\Omega}|\nabla q|$ as $\epsilon\rightarrow0$. In other words, instead of considering the TV model (\ref{eq:ROTV}), the following approximated model with a smoothing regularization term is popularly considered:
\begin{equation}\label{eq:dBMOP_2}
  \begin{aligned}
    \min_{q,u}&~~\Big\{\frac{1}{2}\int_\Omega q|\nabla u-\nabla u_{\delta}|^2\mathrm{d}x+\alpha TV_{\epsilon}(q)\Big\},\\
    \text{s.t.}~&~~-\nabla\cdot(q\nabla u)=f,~~(q,u)\in K\times H_0^1(\Omega).
  \end{aligned}
\end{equation}

To solve \eqref{eq:dBMOP_2} numerically, we choose the following piecewise linear finite element space $V_h$ to discretize the functions $q$ and $u$:
\begin{equation}\label{eq:Vh}
  \begin{aligned}
    V_h&=\big\{v_h\in C(\Omega):v_h|_{\tau}\in\mathcal{P}_1 \text{ for each } \tau\in\mathcal{T}_h\big\},
  \end{aligned}
\end{equation}
where $\mathcal{P}_1$ is the space consisting of polynomials of degrees less than or equal to one, $\mathcal{T}_h$ denotes a regular partition of $\Omega$ into $d$-simplexes, and $h=\max\limits_{\tau\in\mathcal{T}_h}\mathrm{diam}(\tau)$ is the maximal diameter. Then, we obtain a discretized version of \eqref{eq:dBMOP_2} with finite element approximation as
\begin{equation}\label{eq:dBMOP_3}
  \begin{aligned}
    \min_{q_h,u_h}&~~\Big\{\frac{1}{2}\int_\Omega q_h|\nabla u_h-\nabla u_{\delta}|^2\mathrm{d}x+\alpha TV_{\epsilon}(q_h)\Big\},\\
    \text{s.t.}~&~~\big(q_h\nabla u_h,\nabla \phi_h\big)=(f,\phi_h),~\forall~\phi_h\in \bar{V}_h;~~(q_h,u_h)\in K_h\times \bar{V}_h,
  \end{aligned}
\end{equation}
where $K_h:=V_h\cap K$, $\bar{V}_h:=V_h\cap H_0^1(\Omega)$ and $(\cdot,\cdot)$ is the regular inner-product in $L^2(\Omega)$.
For solving the model \eqref{eq:dBMOP_3}, a particularly useful approach is the augmented Lagrangian technique developed in \cite{ItoK:1990a} and then widely used in other literatures such as \cite{ChanT:1997,ChanT:2003,ChenZ:1999,ItoKK:1991,ItoK:1996,KeungZ:1998,TaiFEMC:1998}. More precisely, the augmented Lagrangian functional of \eqref{eq:dBMOP_3} is
\begin{equation}\label{eq:dLagrfunc}
  \begin{aligned}
    {L}_{\gamma_k}(q_h,u_h;\mu_h):&=~\frac{1}{2}\int_\Omega q_h|\nabla u_h-\nabla u_{\delta}|^2\mathrm{d}x+\alpha TV_{\epsilon}(q_h)+(\nabla e(q_h,u_h),\nabla\mu_h)+\frac{\gamma_k}{2}\|\nabla e(q_h,u_h)\|_{L^2(\Omega)}^2,
  \end{aligned}
\end{equation}
with $\gamma_k>0$ the penalty parameter and $\mu_h\in V_h$ the Lagrange multiplier.
The constraint in \eqref{eq:dBMOP_3} is augmented by
\begin{equation}\label{eq:eqv}
  (\nabla e(q_h,u_h),\nabla \phi_h):=(q_h\nabla u_h, \nabla \phi_h)-(f,\phi_h), ~\forall~(q_h,u_h)\in K_h\times \bar{V}_h, ~\forall~\phi_h\in \bar{V}_h,
\end{equation}
where $e(q,u)=(-\Delta)^{-1}\big(\nabla\cdot(q\nabla u)+f\big)$ and $e(\cdot, \cdot)$ can be viewed as an operator from $K\times H_0^1(\Omega)$ to $H_0^1(\Omega)$.
In \cite{ChenZ:1999}, it is proved that the discretized augmented Lagrangian functional \eqref{eq:dLagrfunc} exists at least one saddle-point, and the finite element solution $(u_h,p_h)$ converges to the solution of \eqref{eq:ROTV}. In \cite{ChenZ:1999,ItoK:1990a}, it is suggested to apply the classic augmented Lagrangian method (ALM) originally proposed in \cite{Hestenes:1969,Powell:1969} to (\ref{eq:dBMOP_3}), and the iterative scheme reads as
\begin{equation}\label{eq:ALM}
  \left\{
  \begin{aligned}
    &(q_h^{k+1},u_h^{k+1})=\underset{(q_h,u_h)\in K_h\times\bar{V}_h}{\arg\min}~  {L}_{\gamma_k}(q_h,u_h;\mu_h^k),\\
    &\mu_h^{k+1}=\mu_h^k+\gamma_k e(q_h^{k+1},u_h^{k+1}).\\
  \end{aligned}\right.
\end{equation}
Note that the $(q_h,u_h)$-subproblem in \eqref{eq:ALM} is a smooth optimization problem with a box constraint $K_h$ on the variable $q_h$. In \cite{ChenZ:1999}, global convergence of \eqref{eq:ALM} is proved under the condition that the exact solution of the $(q_h,u_h)$-subproblem of \eqref{eq:ALM} can be obtained at each iteration.

\subsection{Motivations and goals}
Smoothing the TV term loses the originally favorable nonsmoothness property, but enables the eligibility of applying the well known ALM (\ref{eq:ALM}). Meanwhile, the ALM (\ref{eq:ALM}) is mainly of conceptual sense because it is very challenging to implement it numerically. As remarked in \cite{KeungZ:2001}, the smoothing surrogate $TV_{\epsilon}(q)$ leads to a nearly singular and indefinite nonlinear minimization system and solving this system is ``a big difficulty to the numerical resolution process''. Indeed, augmenting the constraint $e(q_h,u_h)=0$ makes the augmented Lagrangian functional ${L}_{\gamma}(q_h,u_h;\mu_h^k)$ nonconvex, and hence the hidden convexity with respect to $q$ in (\ref{eq:ROTV}) is also lost in the ALM (\ref{eq:ALM}). The nonconvex $(q_h,u_h)$-subproblem in the ALM (\ref{eq:ALM}) is numerically difficult also because of the high dimensionality of its variables, the coupling of different variables, as well as its nonlinear structure and ill-conditionedness. In literatures such as \cite{ChenZ:1999,KeungZ:2001}, it is suggested to solve the $(q_h, u_h)$-subproblem inexactly by splitting the variables $q_h$ and $u_h$, and then solving them alternatively. As analyzed in \cite{KeungZ:2001}, the resulting $u_h$-subproblem is a linear yet ill-conditioned saddle-point system and the $q_h$-subproblem is a nearly singular nonlinear minimization problem --- both are still very difficult. It is suggested in \cite{ChenZ:1999} to apply some first-order algorithm with an Armijo line search to solve the decomposed $q_h$-subproblem, each iteration of which also requires solving an ill-conditioned linear system. All these strategies are targeted for approximating the solution of the $(q_h, u_h)$-subproblem in (\ref{eq:ALM}) heuristically, without any guarantee to the theoretically rigorous convergence. All these difficulties become much severer if a higher-dimensional space with $d\ge2$ is considered and fine mesh discretization is used. Indeed, the dimension and condition numbers of the involved linear systems are both of order $O(h^{-d})$. Hence, it is easy to understand the lack of numerical study in the literatures for higher-dimensional spaces of $d\ge 2$ in (\ref{eq:ROTV}) with fine mesh discretization. To the best of our knowledge, only some limited numerical studies for the case where $d=1$ in \eqref{eq:dBMOP_2} and coarse mesh discretization (e.g., $h=1/80$) are available in \cite{ChenZ:1999, Zou:1998}. To summarize, it is extremely challenging to find the exact solution, or even an approximate solution with good accuracy, of the nonconvex, nonlinear, ill-conditioned and large-scaled $(q_h,u_h)$-subproblem in (\ref{eq:ALM}). This challenge posts substantial difficulties to validate the condition in \cite{ChenZ:1999} to guarantee the convergence of the ALM (\ref{eq:ALM}).

Because of the mentioned difficulties in the smoothing surrogate $TV_{\epsilon}(q)$ and the ALM (\ref{eq:ALM}), we are motivated to turn to consider solving the original model (\ref{eq:ROTV}) directly. Our goals are: (1) to tackle the original TV model (\ref{eq:ROTV}) so that the nonsmoothness properties of the diffusion coefficient of (\ref{eq:IDCEE}) can be inherited throughout; (2) to keep the convexity of the data-fidelity functional in \eqref{eq:ROTV} with respect to $q$ throughout; (3) to design an implementable algorithm without any difficult subproblem such as the $(q_h, u_h)$-subproblem in (\ref{eq:ALM}) while it is efficient even for higher-dimensional space of (\ref{eq:ROTV}) with $d=2$ and fine mesh discretization. We will show that the first two goals can be fully achieved by applying the well-known alternating direction method of multipliers (ADMM) which was proposed originally in \cite{GlowinskiM:1975}. For the third goal, we should meticulously investigate the resulting subproblems, and propose some structure-exploiting strategies to tackle these subproblems more effectively. It is mentionable that the curse of dimensionality really matters from the numerical point of view. For example,
for the case where the uniform mesh size $h=1/1024$, the order of dimensionality of the resulting linear systems increases from $10^3$ to $10^6$ if the domain is changed from the unit interval $\Omega \subset \mathbb{R}$ to the unit square $\Omega \subset \mathbb{R}^2$,
while if the domain is fixed as the unit square $\Omega \subset \mathbb{R}^2$, then the order of dimensionality of the resulting linear systems increases from $10^4$ to $10^6$ if the mesh size is refined from $1/128$ to $1/1024$.

\subsection{Conceptual application of ADMM to the original TV model}

As mentioned in \cite{HaoQ:2011}, the elliptic equation \eqref{eq:IDCEE} has a unique weak solution $u$ in $H_0^1(\Omega)$ for each $q\in K$ and $u$ is nonlinearly dependent on $q$. Then, the nonlinear coefficient-to-solution mapping $U:K\rightarrow H_0^1(\Omega)$, which maps each $q\in K$ to the unique solution $u=U(q)\in H_0^1(\Omega)$ of \eqref{eq:IDCEE}, is well defined. Instead of augmenting the elliptic equation \eqref{eq:IDCEE} as a constraint by introducing $e(q,u)$ in \eqref{eq:eqv}, we temporarily take the liberty to represent $u$ as a function of $q$ via the equation \eqref{eq:IDCEE}, denote by
$$
J(q)=:\frac{1}{2}\int_\Omega q|\nabla U(q)-\nabla u_{\delta}|^2\mathrm{d}x,
$$
and then reformulate the model \eqref{eq:ROTV} as a minimization problem only depending on $q$. Then the finite element discretized version of the nonsmooth problem \eqref{eq:ROTV} can be written as
\begin{equation}\label{eq:dROTV}
  \min_{q_h\in K_h}\Big\{J(q_h)+\alpha\|\nabla q_h\|_{L^1(\Omega)}\Big\}.
\end{equation}

To implement the ADMM to solve (\ref{eq:dROTV}), there are multiple ways. For instance, it is easy to consider introducing an auxiliary variable $p_h:=\nabla q_h$ so as to replace $\nabla q_h$ with $p_h$ in the objective functional of \eqref{eq:dROTV}. That is, the model \eqref{eq:dROTV} can be reformulated as
\begin{equation*}
  \begin{aligned}
    \min_{q_h,p_h}&~\Big\{J(q_h)+\alpha\|p_h\|_{L^1(\Omega)}\Big\},\\
    \text{s.t.}~&~~\nabla q_h-p_h=0,~~(q_h,p_h)\in K_h\times W_h,
  \end{aligned}
\end{equation*}
where $W_h:=\{p_h\in L^1(\Omega;\mathbb{R}^d):p_h|_{\tau}~ \textrm{is constant for each } \tau\in \mathcal{T}_h\}$.
Inspired by \cite{BartelsM:2016}, we can employ the weighted $L^2$-inner product $(\cdot,\cdot)_h=h^d(\cdot,\cdot)$ and the corresponding norm $\|\cdot\|_h=h^{d/2}\|\cdot\|_{L^2(\Omega)}$ with $d$ the space dimension, to penalize the constraint. Then, the augmented Lagrangian functional is
\begin{equation}\label{eq:AL-new}
  \begin{aligned}
    \hat{L}_{\beta}(q_h,p_h;\lambda_h):=~J(q_h)+\alpha\|p_h\|_{L^1(\Omega)}
    +(\nabla q_h-p_h,\lambda_h)_h+\frac{\beta}{2}\|\nabla q_h-p_h\|_h^2,
  \end{aligned}
\end{equation}
where $\beta>0$ is a penalty parameter. The corresponding ADMM iterative scheme reads as
\begin{equation}\label{eq:ADMM_b}
  \left\{
  \begin{aligned}
    &q_h^{k+1}=\underset{q_h\in K_h}{\arg\min}~\hat{L}_{\beta}(q_h,p_h^k;\lambda_h^k),\\
    &p_h^{k+1}=\underset{p_h\in W_h}{\arg\min}~\hat{L}_{\beta}(q_h^{k+1},p_h;\lambda_h^k),\\
    &\lambda_h^{k+1}=\lambda_h^k+\beta(\nabla q_h^{k+1}-p_h^{k+1}).\\
  \end{aligned}\right.
\end{equation}
In \eqref{eq:AL-new}, we do not use the regular $L^2$-inner product $(\cdot, \cdot)$ and its induced $L^2$-norm penalty term $\frac{\beta}{2}\|\nabla q_h-p_h\|_{L^2(\Omega)}^2$. Indeed, as analyzed in \cite{BartelsM:2016}, the $L^2$-inner product may lead to numerical instability because $\|\nabla q_h\|_{L^2(\Omega)}$ is unbounded. It is further noticed in \cite{BartelsM:2016} that an inverse estimate shows that $\nabla q_h$ is bounded with respect to $\|\cdot\|_h$ and the corresponding scheme tends to be more numerically stable.

Note that the gradient operator $\nabla$ is involved in the penalty term in \eqref{eq:AL-new}, and as analyzed in \cite{LarsonB:2013}, the condition number of the corresponding stiffness matrix (whose entries are $(\nabla\phi_h^j,\nabla\phi_h^i)$ with $\phi_h^i$ the finite element basis functions of $V_h$) is of order $h^{-d}$. Hence, the condition number of the stiffness matrix may be extremely high for fine mesh. Accordingly, numerical performance of (\ref{eq:ADMM_b}) may be more severely affected by the penalty parameter $\beta$ if fine mesh is used for discretization, which can be easily verified by numerical experiments. Because of this concern, we prefer to penalize some term irrelevant to the gradient operator $\nabla$. Note that the condition number of the mass matrix is bounded and independent of the mesh size $h$; see, e.g., \cite{LarsonB:2013}. Therefore, we introduce the auxiliary variable $p_h$ and replace $q_h$ in the TV term. That is, we reformulate the TV model \eqref{eq:dROTV} as
\begin{equation}\label{eq:dROTV2}
  \begin{aligned}
    \min_{q_h,p_h}&~\Big\{J(q_h)+\alpha\|\nabla p_h\|_{L^1(\Omega)}\Big\},\\
    \text{s.t.}~&~~q_h-p_h=0,~~(q_h,p_h)\in K_h\times V_h.
  \end{aligned}
\end{equation}
The corresponding augmented Lagrangian functional of \eqref{eq:dROTV2} is
\begin{equation}\label{eq:dLagrfunc_2}
  \begin{aligned}
    {\mathcal{L}}_\beta(q_h,p_h;\lambda_h):=~J(q_h)+\alpha\| \nabla p_h\|_{L^1(\Omega)}+( q_h-p_h,\lambda_h)+\frac{\beta}{2}\| q_h-p_h\|_{L^2(\Omega)}^2,
  \end{aligned}
\end{equation}
and the corresponding ADMM scheme reads as:
\begin{subequations}\label{eq:ADMM_a_origin}
  \begin{align}
    &q_h^{k+1}=\underset{q_h\in K_h}{\arg\min}~ {\mathcal{L}}_\beta(q_h,p_h^k;\lambda_h^k),\tag{\ref{eq:ADMM_a_origin}a}\label{eq:ADMM_a_origin_a}\\
    &p_h^{k+1}=\underset{p_h\in V_h}{\arg\min}~ {\mathcal{L}}_\beta(q_h^{k+1},p_h;\lambda_h^k),\tag{\ref{eq:ADMM_a_origin}b}\label{eq:ADMM_a_origin_b}\\
    &\lambda_h^{k+1}=\lambda_h^k+\beta(q_h^{k+1}-p_h^{k+1}).\tag{\ref{eq:ADMM_a_origin}c}
  \end{align}\vspace*{-90pt}
  \begin{equation*}\hspace*{-160pt}
    \left\{
    \begin{aligned}
      \\ \\ \\ \\[5pt]
    \end{aligned}
    \right.
  \end{equation*}
\end{subequations}

It is arguably trivial to derive the ADMM \eqref{eq:ADMM_a_origin} conceptually. But it is clear that both the subproblems in (\ref{eq:ADMM_a_origin}) are convex and the scheme is for solving the discretized version of the original TV model \eqref{eq:ROTV}. Hence, the ADMM (\ref{eq:ADMM_a_origin}) essentially differs from the ALM \eqref{eq:ALM} in the sense that the original TV term as well as the convexity with respect to the variable $q_h$ are both kept. On the other hand, as we shall show in Section \ref{sec:ADMM}, despite its convexity, it is highly nontrivial to solve the resulting subproblems, especially the $q_h$-subproblem (\ref{eq:ADMM_a_origin_a}). Hence, the ADMM \eqref{eq:ADMM_a_origin} is numerically meaningful only if both the $q_h$- and $p_h$-subproblems can be solved efficiently, especially for the case where $d\ge2$ and $h$ is small.

\subsection{Organization}

The rest of this paper is organized as follows. In Section \ref{sec:pre}, some preliminaries are summarized for further analysis. Then, we focus on the subproblems (\ref{eq:ADMM_a_origin_a}) and (\ref{eq:ADMM_a_origin_b}) in Sections \ref{sec:ADMM} and \ref{sec:CNN}, respectively. The flowchart of implementation of the proposed numerical approach is presented in Section \ref{sec:ADMM-Newton-CNN}. Some preliminary numerical results are reported in Section \ref{sec:num} to verify the efficiency of the proposed numerical approach. Finally, some conclusions are drawn in Section \ref{sec:Conc}.

\section{Preliminaries}\label{sec:pre}
In this section, we summarize some preliminaries which will be used for further analysis.
We say that $u\in H_0^1(\Omega)$ is a weak solution of the elliptic equation \eqref{eq:IDCEE} if it satisfies the following variational form:
$$
  a(u,\phi):=\int_{\Omega}q\nabla u\cdot\nabla\phi\mathrm{d}x=\int_{\Omega}f\phi\mathrm{d}x,~~\forall~\phi\in H_0^1(\Omega).
$$
The bilinear form $a(\cdot,\cdot)$ satisfies the coercivity condition $a(u,u)\ge C_{\Omega}\|u\|_{H_0^1(\Omega)}^2$ for any $u\in H_0^1(\Omega)$ and $q\in K$, where $C_{\Omega}$ is a positive constant depending on $\Omega$ and the low bound $a_0$ of $q$. By the Lax-Milgram theorem \cite{Evans:2010}, there exists a unique weak solution $u$ of \eqref{eq:IDCEE} in $H_0^1(\Omega)$ for each $q\in K$, and $u$ is nonlinearly dependent on $q$. Then, we can define the nonlinear coefficient-to-solution mapping $U:K\rightarrow H_0^1(\Omega)$, which maps each $q\in K$ to the unique solution $u=U(q)\in H_0^1(\Omega)$ of \eqref{eq:IDCEE}; see \cite{HaoQ:2011} for more details.

For discretization, because of the low regularity of the functions in the space $BV(\Omega)$, only some low order polynomials will be chosen for the finite element space. As mentioned in \cite{Bartels:2012}, the piecewise affine globally continuous finite element spaces are dense in $BV(\Omega)$ with respect to weak* convergence in $BV(\Omega)$, while in general the piecewise constant finite element approximation for $u$ cannot be expected to converge to an exact solution. Thus, we discretize the model \eqref{eq:ROTV} in the finite element space $V_h$ (see \eqref{eq:Vh}) and obtain the following discrete problem:
\begin{equation}\label{eq:dBMOP}
  \min_{q_h\in K_h}\Big\{J(q_h)+\alpha\|\nabla q_h\|_{L^1(\Omega)}\Big\},
\end{equation}
where $K_h=V_h\cap K$, and $V_h$ and $K$ are given by \eqref{eq:Vh} and \eqref{eq:K}, respectively. The energy functional $J(q_h)$ of \eqref{eq:dBMOP} is
\begin{equation}\label{eq:dBMOP_q}
  J(q_h):=\frac{1}{2}\int_{\Omega}q_h|\nabla U(q_h)- \nabla u_{\delta} |^2\mathrm{d}x,
\end{equation}
where $U(q_h)$ is the solution of the following variational form:
\begin{equation*}
  a_h(q_h,U(q_h);\phi_h)=( f,\phi_h),~\forall~\phi_h\in \bar{V}_h,
\end{equation*}
with $\bar{V}_h:=V_h\cap H_0^1(\Omega)$ and
\begin{equation}
  a_h(q_h,U(q_h);\phi_h):=\int_{\Omega}q_h\nabla U(q_h)\cdot\nabla\phi_h \mathrm{d}x.
\end{equation}

According to \cite[Lemma 2.3]{HaoQ:2010}, the functional $J(\cdot)$ in \eqref{eq:dBMOP_q} is convex on the convex set $K$. For any $q\in K, \xi, \eta\in L^\infty(\Omega)$, the first derivative of $J(\cdot)$ is given by
\begin{equation}\label{eq:FDJa}
  J'(q)\xi=-\frac{1}{2}\int_{\Omega} \xi|\nabla U(q)|^2 \mathrm{d}x+\frac{1}{2}\int_{\Omega} \xi |\nabla u_{\delta}|^2 \mathrm{d}x,
\end{equation}
and the second derivative of $J(\cdot)$ is given by
\begin{equation}\label{eq:SDJa}
  J''(q)(\xi,\eta)=-\int_{\Omega} \xi\nabla U(q)\cdot \nabla U'(q) \eta \mathrm{d}x,
\end{equation}
where $U'(q) \eta$ satisfies
\begin{equation}\label{eq:U1Ja}
  \int_{\Omega} q \nabla U'(q)\eta\cdot\nabla v\mathrm{d}x=- \int_{\Omega} \eta\nabla U(q)\cdot\nabla v\mathrm{d}x,~\forall~v\in H^1_0(\Omega).
\end{equation}

In the next two sections, we will elaborate on how to solve the resulting subproblems for the ADMM \eqref{eq:ADMM_a_origin}.
For notational convenience, we denote by  $\bm{q}:=(\bm{q}_1,\dots, \bm{q}_N)^\top\in \mathbb{R}^N$ the coefficients of $q_h$. That is, $q_h=\sum_{i=1}^N \bm{q}_i \phi_h^i$, where $\{\phi_h^i\}_{i=1}^N$ are the finite element basis functions in $V_h$. The same setting is also applied to $p_h,~\lambda_h,~u_h~,r_h$, with the coefficients $\bm{p},~\bm{\lambda},~\bm{u},~\bm{r}\in  \mathbb{R}^N$, respectively. Then, we define the function $J(\bm{q}):=J(q_h)$.

\section{Active-set Newton method for the $q_h$-subproblem}\label{sec:ADMM}

In this section, we focus on the $q_h$-subproblem \eqref{eq:ADMM_a_origin_a}. How to solve this subproblem is crucial to ensure the performance of the ADMM (\ref{eq:ADMM_a_origin}), and this is the most technical part of the paper.

\subsection{Optimality conditions}

For the $q_h$-subproblem \eqref{eq:ADMM_a_origin_a}, it can be written as the following smooth and nonlinear optimization problem:
\begin{equation}\label{eq:qh}
  q_h^{k+1}=\underset{q_h\in K_h}{\arg\min}~ \Big\{J(q_h)+\frac{\beta}{2}\|q_h-p_h^k+\frac{\lambda_h^k}{\beta}\|_{L^2(\Omega)}^2\Big\}.
\end{equation}
We further reformulate the problem \eqref{eq:qh} in Euclidean space and derive its first-order optimality conditions.
Recall the definitions of $\bm{q},~\bm{p},~\bm{\lambda}$ and $J(\bm{q})$ in Section \ref{sec:pre}.
The optimization problem \eqref{eq:qh} can be rewritten as
\begin{equation}\label{qh_dis}
  \begin{aligned}
    &\min_{\bm{q}\in \mathbb{R}^N}\Big\{J(\bm{q})+\frac{\beta}{2}\| \bm{q}-\bm{p}^k+\frac{\bm{\lambda}^k}{\beta}\|_{M}^2\Big\},\\
    &~~~\text{s.t.}\quad \bm{a}_0 \le \bm{q}\le \bm{a}_1,
  \end{aligned}
\end{equation}
where $M$ denotes the mass matrix as $M_{i,j} = (\phi_h^j,\phi_h^i)$, $\bm{a}_0=a_0\bm{1}$, and $\bm{a}_1=a_1\bm{1}$. Here, $\bm{1}\in\mathbb{R}^N$ denotes the vector with constant entries $1$.
The Lagrangian function of (\ref{qh_dis}) is
$$
  l(\bm{q},\bm{\eta}_0,\bm{\eta}_1)=J(\bm{q})+\frac{\beta}{2}\| \bm{q}-\bm{p}^k+\frac{\bm{\lambda}^k}{\beta}\|_{M}^2+\bm{\eta}_0^\top(\bm{q}-\bm{a}_0)+\bm{\eta}_1^\top(\bm{q}-\bm{a}_1),
$$
with $\bm{\eta}_0, \bm{\eta}_1\in \mathbb{R}^N$ the Lagrange multipliers.
Then the corresponding KKT conditions are
$$\left \{
  \begin{aligned}
	&J'(\bm{q})+\beta M(\bm{q}-\bm{p}^k+\frac{\bm{\lambda}^k}{\beta})+(\bm{\eta}_0+\bm{\eta}_1)=0,\\
    &\bm{\eta}_1\ge0,~ \bm{q}\le \bm{a}_1,~ \bm{\eta}_1^\top(\bm{q}-\bm{a}_1)=0,\\
    & \bm{\eta}_0\le0,~ \bm{q}\ge \bm{a}_0,~ \bm{\eta}_0^\top(\bm{q}-\bm{a}_0)=0.\\
  \end{aligned}\right.
$$
Furthermore, denoting $\bm{\eta}:=\bm{\eta}_0+\bm{\eta}_1$ and
$$
  C( \bm{q}, \bm{\eta}):=\bm{\eta}-\max \{0, \bm{\eta}+c(\bm{q}-\bm{a}_1)\}-\min\{0, \bm{\eta}+c(\bm{q}-\bm{a}_0)\}, ~c>0,
$$
we can represent the KKT conditions as the equation
\begin{equation}\label{eq:KKT}
  F(\bm{q},\bm{\eta}):=
  \begin{pmatrix}
    J'(\bm{q})+\beta M(\bm{q}-\bm{p}^k+\frac{\bm{\lambda}^k}{\beta})+\bm{\eta} \\
    C(\bm{q}, \bm{\eta})
  \end{pmatrix}=0.
\end{equation}

\subsection{Computation of the first-order derivative}\label{subS_gradient}

To solve the $q_h$-subproblem (\ref{eq:qh}), it is natural to consider the first-order derivative of $J(q_h)$ and probe its computational complexity.
It follows from \eqref{eq:FDJa} that the first-order derivative of $J(q_h)$ satisfies
$$
  (\xi_h, J'(q_h))=(\xi_h, -\frac{1}{2} |\nabla U(q_h)|^2  +\frac{1}{2}|\nabla u_{\delta}|^2), ~\forall~\xi_h \in V_h,
$$
where $U(q_h)$ is the solution to
\begin{equation}\label{u_subproblem}
  a_h(q_h,U(q_h);\phi_h)=(f,\phi_h),~\forall~\phi_h\in \bar{V}_h=V_h\cap H_0^1(\Omega).
\end{equation}
Let $k_m$ with $m\ge 0$ be iteration counter for the inner loop for solving the $q_h$-subproblem at the $k$-th iteration; $k_0$ be the initial iterate for the inner loop. Then, how to compute ${J'}(\bm{q}^{k_m})$ can be summarized in the following Subroutine \ref{Compute_J'a}.

\begin{algorithm}[htpb]\caption{Computation of ${J'}(\bm{q}^{k_m})$.}\label{Compute_J'a}
  \begin{algorithmic}[1]
    \Function {Gradient}{$\bm{q}^{k_m}$}
      \State Obtain $\bar{\bm{u}}^{k_m}$ via solving
        \begin{equation}\label{eq:A_{k_m}}
          A_{k_m} \bar{\bm{u}}^{k_m}=\bm{f},
        \end{equation}
        where $(A_{k_m})_{i,j}=(q_h^{k_m} \nabla\phi_h^j, \nabla\phi_h^i)$ and $\bm{f}$ defined as $\bm{f}_i=(f,\phi_h^i)$.
      \State Substitute $U(q_h^{k_m})=\sum_{j=1}^N \bar{\bm{u}}_j^{k_m}\phi_h^j$ into the following equation to compute ${J'}(\bm{q}^{k_m})$
        $$
        \begin{aligned}
          ({J'}(\bm{q}^{k_m}))_i&=(\phi_h^i,  J'(q_h^{k_m}))=( \phi_h^i, -\frac{1}{2} |\nabla U(q_h^{k_m})|^2  +\frac{1}{2}|\nabla u_{\delta}|^2), ~ i=1,\dots, N.
        \end{aligned}$$
      \State \Return{${J'}(\bm{q}^{k_m})$.}
	\EndFunction			
  \end{algorithmic}
\end{algorithm}

It is easy to see that the computation of $J'(q_h)$ requires values of $U(q_h)$, which should be obtained by computing $A_{k_m}$ and solving the linear system \eqref{eq:A_{k_m}} iteratively. Note that the linear system \eqref{eq:A_{k_m}} is a discretized formulation of the elliptic equation \eqref{eq:IDCEE}. As analyzed in \cite{LarsonB:2013}, its dimension and the condition number of the coefficient matrix $A_{k_m}$ are both of order $O(h^{-d})$. Thus the linear system \eqref{eq:A_{k_m}} is large-scaled and ill-conditioned for discretization with fine mesh, and computing $J'(q_h)$ may be expensive. Note that computing the objective function value in \eqref{eq:qh} requires values of $U(q_h)$ as well. Hence, these difficulties essentially imply that it is computationally demanding even if some first-order algorithm is applied to seek a medium- or low-accuracy numerical solution of the problem (\ref{eq:qh}). Indeed, implementing a first-order algorithm usually requires certain line-search techniques with multiple computations of the objective function values, to discern appropriate step sizes. Our numerical experiments actually validate the failure of a number of popular first-order algorithms (such as the gradient projection method and the conjugate gradient projection method with backtracking line-search) firmly for solving the subproblem (\ref{eq:qh}).

\subsection{Active-set Newton method for the problem \eqref{eq:qh}}\label{ssec:acn}

As analyzed, though the ADMM (\ref{eq:ADMM_a_origin}) per se can be easily derived, it is keen to solve the $q_h$-subproblem (\ref{eq:ADMM_a_origin_a}), i.e., the problem \eqref{eq:qh}. Because demanding computation is required yet only a medium- or low-accuracy solution can be targeted, it is not attractive to consider first-order algorithms for this subproblem. It is thus interesting to investigate how much more complicated if a second-order algorithm is applied to the problem \eqref{eq:qh}. In this and the next subsections, we will show that, counter-intuitively, the Newtonian system of \eqref{eq:KKT} can be appropriately reformulated and relaxed so that its computation reduces to solving a simple positive definite linear system, and then the benchmark active-set Newton method in, e.g., \cite{HinzeP:2008,KR2002}, can be applied very efficiently. Computation of the Newton step is comparable with, and usually less than, that of a single iteration of the gradient projection method with some backtracking line-search strategy, while the accuracy is much higher. This is a convincing example of deriving model-tailored efficient algorithms by taking full advantage of the structure of the model under discussion.

To elaborate on the active-set Newton method for \eqref{eq:KKT}, let us define
$$
  {\mathcal A}_{k_m}={\mathcal A}^+_{k_m} \cup {\mathcal A}^-_{k_m} \quad \text{and}
  \quad {\mathcal I}_{k_m}=\{1,\dots, N\}\setminus {\mathcal A}_{k_m}
$$
as the sets of the active and inactive indices at $(\bm{q}^{k_m},\bm{\eta}^{k_m})$, respectively, where ${\mathcal A}^+_{k_m}$ and ${\mathcal A}^-_{k_m}$ are the sets given respectively by
$$
  {\mathcal A}_{k_m}^+=\{i~|\bm{\eta}_i^{k_m}+c(\bm{q}^{k_m}-\bm{a}_1)_i>0\}\quad \text{and}
  \quad {\mathcal A}_{k_m}^-=\{i~|\bm{\eta}_i^{k_m}+c(\bm{q}^{k_m}-\bm{a}_0)_i<0\}.
$$
For the mapping $F(\bm{q}^{k_m}, \bm{\eta}^{k_m})$ defined in \eqref{eq:KKT}, let $\partial F(\bm{q}^{k_m}, \bm{\eta}^{k_m})$ be the generalized Jacobian of \eqref{eq:KKT} in sense of Clarke (see \cite{Clarke:1983}). Then, as analyzed in \cite{KR2002}, we have
\begin{equation}\label{partial_F}
 F'(\bm{q}^{k_m}, \bm{\eta}^{k_m}):= \begin{pmatrix}
    J''(\bm{q}^{k_m})+\beta M&I\\
    -c\Pi_{{\mathcal A}_{k_m}}& \Pi_{{\mathcal I}_{k_m}}
  \end{pmatrix} \in \partial F(\bm{q}^{k_m}, \bm{\eta}^{k_m}),
\end{equation}
where $\Pi_{{\mathcal A}_{k_m}}$ and $\Pi_{{\mathcal I}_{k_m}}$ denote the diagonal binary matrices with nonzero entries in ${\mathcal A}_{k_m}$ and ${\mathcal I}_{k_m}$, respectively.
With \eqref{partial_F}, it is easy to see that the $k_m$-th iteration of the active-set Newton method for \eqref{eq:KKT} is solving
\begin{equation}\label{Full_Newton}\small
  \begin{pmatrix}
    J''(\bm{q}^{k_m})+\beta M&I\\
    -c\Pi_{{\mathcal A}_{k_m}}& \Pi_{{\mathcal I}_{k_m}}
  \end{pmatrix}
  \begin{pmatrix}
    \bm{q}^{k_{m+1}}-  \bm{q}^{k_{m}}\\\bm{\eta}^{k_{m+1}}-\bm{\eta}^{k_{m}}
  \end{pmatrix}=-
  \begin{pmatrix}
{J'}(\bm{q}^{k_m})+\beta M(\bm{q}^k-\bm{p}^k+\frac{\bm{\lambda}^k}{\beta})\\ \bm{\eta}^{k_{m}}- \Pi_{{\mathcal A}_{k_m}^+}
   (\bm{\eta}^{k_{m}}+c(\bm{q}^{k_m}-\bm{a}_1))-\Pi_{{\mathcal A}_{k_m}^-}(\bm{\eta}^{k_{m}}+c(\bm{q}^{k_m}-\bm{a}_0))
  \end{pmatrix}.
\end{equation}
Since the second equation of \eqref{Full_Newton} implies that
$$
(\bm{\eta}^{k_{m+1}})_{{\mathcal I}_{k_m}}=(\bm{\eta}^{k_{m+1}}-\bm{\eta}^{k_{m}})_{{\mathcal I}_{k_m}}+(\bm{\eta}^{k_{m}})_{{\mathcal I}_{k_m}}=0,
$$
we can remove those rows that belong to the indices in ${\mathcal I}_{k_m}$ from the second equation of \eqref{Full_Newton}, and simplify \eqref{Full_Newton} as
\begin{equation}\label{eq:Newton_original}
  \begin{pmatrix}
    \beta  J''(\bm{q}^{k_m})+\beta M&{\mathcal P}_{{\mathcal A}_{k_m}}^\top\\
    {\mathcal P}_{{\mathcal A}_{k_m}}&0
  \end{pmatrix}
  \begin{pmatrix}
    \bm{q}^{k_{m+1}}-\bm{q}^{k_m}\\(\bm{\eta}^{k_{m+1}})_{{\mathcal A}_{k_m}}
  \end{pmatrix}=
  \begin{pmatrix}
    -{J'}(\bm{q}^{k_m})-\beta M(\bm{q}^{k_m}-\bm{p}^k+\frac{\bm{\lambda}^k}{\beta})\\{\mathcal P}_{{\mathcal A}_{k_m}^+}\bm{a}_1+{\mathcal P}_{{\mathcal A}_{k_m}^-}\bm{a}_0-{\mathcal P}_{{\mathcal A}_{k_m}}\bm{q}^{k_m}
  \end{pmatrix}.
\end{equation}
In (\ref{eq:Newton_original}), ${\mathcal P}_{\mathcal C}$ denotes the matrix consisting of those rows of $\Pi_{\mathcal C}$ that belong to the indices in a given set ${\mathcal C}$. Obviously, it holds that $\Pi_{\mathcal C}={\mathcal P}_{\mathcal C}^\top {\mathcal P}_{\mathcal C}$.

Though it is trivial to analytically derive the system \eqref {eq:Newton_original} for the $k_m$-th iteration of the active-set Newton method for \eqref{eq:KKT}, how to solve  \eqref{eq:Newton_original} numerically deserves meticulous analysis mainly because computing the second-order derivative $J''(\bm{q}^{k_m})$ is very expensive. Indeed, it follows from \eqref{eq:SDJa} that computing the second-order derivative $J''(\bm{q}^{k_m})$ directly at each iteration requires computing $\{U'(q_h^{k_m}) \phi_h^i\}_{i=1}^N$ from \eqref{eq:U1Ja} for each finite element basis function $\phi_h^i$. This means a sequence of discretized elliptic equations in form of \eqref{eq:A_{k_m}} are needed to be solved, and recall that each of them is large-scaled and ill-conditioned for fine mesh discretization.

To avoid computing $J''(\bm{q}^{k_m})$, we take an alternative approach to compute $\bm{q}^{k_{m+1}}-\bm{q}^{k_m}$. The key idea is substituting the discrete equations of \eqref{eq:SDJa} and \eqref{eq:U1Ja} into \eqref{eq:Newton_original} to eliminate $J''(\bm{q}^{k_m})$.
To see the details, it follows from \eqref{eq:SDJa} that
$$
  J''(q_h^{k_m})(\xi_h, q_h^{k_{m+1}}-q_h^{k_m})=-(\xi_h \nabla U(q_h^{k_m}), \nabla U'(q_h^{k_m})(q_h^{k_{m+1}}-q_h^{k_m})).
$$
Then, we have
$$
J''(\bm{q}^{k_m})(\bm{q}^{k_{m+1}}-\bm{q}^{k_m})=-N_{k_m} \bm{r}^{k_{m+1}},
$$
where
$$
(N_{k_m})_{i,j}=(\phi_h^j \nabla U(q_h^{k_m}), \nabla\phi_h^i),
\;\; r_h^{k_{m+1}}= U'(q_h^{k_m})(q_h^{k_{m+1}}-q_h^{k_m})=\sum_{i=1}^N \bm{r}_i^{k_{m+1}} \phi_h^i
$$
and
$$
\bm{r}^{k_{m+1}}=(\bm{r}_1^{k_{m+1}},\dots,\bm{r}_N^{k_{m+1}})^\top.
$$
It also follows from \eqref{eq:U1Ja} that
$$
  a_h(q_h^{k_m},r_h^{k_{m+1}};v_h) =-((q_h^{k_{m+1}}-q_h^{k_m}) \nabla U(q_h^{k_m}), \nabla v_h),~\forall~v_h\in \bar{V}_h,
$$
which implies that
$$
A_{k_m} \bm{r}^{k_{m+1}}=-N_{k_m}^\top(\bm{q}^{k_{m+1}}-\bm{q}^{k_m}).
$$
Thus we have
$$
J''(\bm{q}^{k_m})(\bm{q}^{k_{m+1}}-\bm{q}^{k_m})=-N_{k_m} \bm{r}^{k_{m+1}} \;\hbox{with}\; A_{k_m} \bm{r}^{k_{m+1}}=-N_{k_m}^\top (\bm{q}^{k_{m+1}}-\bm{q}^{k_m}).
$$
Next, substituting
$$
J''(\bm{q}^{k_m})(\bm{q}^{k_{m+1}}-\bm{q}^{k_m})=-N_{k_m} \bm{r}^{k_{m+1}}
$$
into the Newtonian system \eqref{eq:Newton_original}, we obtain the under-determined linear system
\begin{equation}\label{eq:Newton_substituted}
  \begin{pmatrix}
    -\beta N_{k_m} \bm{r}^{k_{m+1}}\\0
  \end{pmatrix}+
  \begin{pmatrix}
    \beta M&{\mathcal P}_{{\mathcal A}_{k_m}}^\top\\
    {\mathcal P}_{{\mathcal A}_{k_m}}&0
  \end{pmatrix}
  \begin{pmatrix}
    \bm{q}^{k_{m+1}}-\bm{q}^{k_m}\\(\bm{\eta}^{k_{m+1}})_{{\mathcal A}_{k_m}}
  \end{pmatrix}=
  \begin{pmatrix}
    -{J'}(\bm{q}^{k_m})-\beta M(\bm{q}^{k_m}-\bm{p}^k+\frac{\bm{\lambda}^k}{\beta})\\{\mathcal P}_{{\mathcal A}_{k_m}^+}\bm{a}_1+{\mathcal P}_{{\mathcal A}_{k_m}^-}\bm{a}_0-{\mathcal P}_{{\mathcal A}_{k_m}}\bm{q}^{k_m}
  \end{pmatrix}
\end{equation}
with respect to $(\bm{r}^{k_{m+1}}, \bm{q}^{k_{m+1}}-\bm{q}^{k_m}, (\bm{\eta}^{k_{m+1}})_{{\mathcal A}_{k_m}})$. Then, combining \eqref{eq:Newton_substituted} with
$$
A_{k_m} \bm{r}^{k_{m+1}}=-N_{k_m}^\top(\bm{q}^{k_{m+1}}-\bm{q}^{k_m}),
$$
we obtain the following expanded linear system:
\begin{equation}\label{Newton_a0}
  \begin{pmatrix}
    A_{k_m} &N_{k_m}^\top&0\\
    -N_{k_m} &\beta M&{\mathcal P}_{{\mathcal A}_{k_m}}^\top\\
    0&{\mathcal P}_{{\mathcal A}_{k_m}}&0
  \end{pmatrix}
  \begin{pmatrix}
    \bm{r}^{k_{m+1}}\\ \bm{q}^{k_{m+1}}-\bm{q}^{k_m}\\(\bm{\eta}^{k_{m+1}})_{{\mathcal A}_{k_m}}
  \end{pmatrix}=
  \begin{pmatrix}0\\
    -{J'}(\bm{q}^{k_m})-\beta M(\bm{q}^{k_m}-\bm{p}^k+\frac{\bm{\lambda}^k}{\beta})\\{\mathcal P}_{{\mathcal A}_{k_m}^+}\bm{a}_1+{\mathcal P}_{{\mathcal A}_{k_m}^-}\bm{a}_0-{\mathcal P}_{{\mathcal A}_{k_m}}\bm{q}^{k_m}
  \end{pmatrix},
\end{equation}
which is equivalent to the Newtonian system \eqref{eq:Newton_original}. Note that there is no need to compute the usually expensive $J''(\bm{q}^{k_m})$ in \eqref{Newton_a0}, and all the matrices $A_{k_m}, N_{k_m}, M, {\mathcal P}_{{\mathcal A}_{k_m}^+}$ and ${\mathcal P}_{{\mathcal A}_{k_m}^-}$ are easy to compute. For convenience, we denote
\begin{equation}\label{F_{k_m}_d_{k_m}}
F^{k_m}:=
  \begin{pmatrix}
    A_{k_m} &N_{k_m}^\top&0\\
    -N_{k_m} &\beta M&{\mathcal P}_{{\mathcal A}_{k_m}}^\top\\
    0&{\mathcal P}_{{\mathcal A}_{k_m}}&0
  \end{pmatrix}{\text{~and~}} \begin{pmatrix}\bm{d}_1\\\bm{d}_2\\\bm{d}_3\end{pmatrix}:=\begin{pmatrix}0\\
  -{J'}(\bm{q}^{k_m})-\beta M(\bm{q}^{k_m}-\bm{p}^k+\frac{\bm{\lambda}^k}{\beta})\\{\mathcal P}_{{\mathcal A}_{k_m}^+}\bm{a}_1+{\mathcal P}_{{\mathcal A}_{k_m}^-}\bm{a}_0-{\mathcal P}_{{\mathcal A}_{k_m}}\bm{q}^{k_m}\end{pmatrix}.
\end{equation}

\subsection{Schur complement reduction}\label{sec:NtsA_II}

Recall that the Newtonian system \eqref{Newton_a0} is an expanded system of the linear saddle-point problem \eqref{eq:Newton_original}, and it is clear that \eqref{Newton_a0} is indefinite. Moreover, because of the stiffness matrix $A_{k_m}$ in its coefficient matrix $F^{k_m}$, the system \eqref{Newton_a0} is also ill-conditioned. Hence, it is not easy to solve the Newtonian system \eqref{Newton_a0}. As analyzed in \cite{BenziG:2005}, there are two types of algorithms that can be used to solve \eqref{Newton_a0}: the \textit{segregated}
and \textit{coupled} (also known as ``all at once") methods.

Note that the right-bottom $2\times 2$ block of the matrix $F^{k_m}$ in \eqref{F_{k_m}_d_{k_m}}
\begin{equation}\label{RB_matrix}
\begin{pmatrix}\beta M&{\mathcal P}_{{\mathcal A}_{k_m}}^\top\\ {\mathcal P}_{{\mathcal A}_{k_m}}&0\end{pmatrix}
\end{equation}
is well-conditioned and hence the variables $(\bm{q}^{k_{m+1}}-\bm{q}^{k_m})$ and $(\bm{\eta}^{k_{m+1}})_{{\mathcal A}_{k_m}}$ can be computed easily once $\bm{r}^{k_{m+1}}$ is obtained.
We are thus inspired to choose the Schur complement reduction in \cite{BenziG:2005}, which is a major segregated approach, to convert the Newtonian system \eqref{Newton_a0} to a linear system with only respect to the variable $\bm{r}^{k_{m+1}}$, by using the block factorization of the coefficient matrix $F^{k_m}$ in \eqref{F_{k_m}_d_{k_m}}. For the matrix $F^{k_m}$, it can be factorized as
\begin{equation}\label{bd_F_{k_m}}
F^{k_m}=
\begin{pmatrix}
I &\frac{1}{\beta}N_{k_m}^\top M^{-1}&R_{k_m}\\
0&I&0\\
0&\frac{1}{\beta}{\mathcal P}_{{\mathcal A}_{k_m}}M^{-1}&I\end{pmatrix}\begin{pmatrix}
S_{k_m}&0&0\\
0 &\beta M&0\\
0&0&-\frac{1}{\beta}({\mathcal P}_{{\mathcal A}_{k_m}}M^{-1}{\mathcal P}_{{\mathcal A}_{k_m}}^\top)\end{pmatrix}\begin{pmatrix}
I &0&\\
-\frac{1}{\beta}M^{-1}N_{k_m}&I&\frac{1}{\beta}M^{-1}{\mathcal {\mathcal P}_{{\mathcal A}_{k_m}}}\\
-R_{k_m}^\top&0&I\end{pmatrix},
\end{equation}
where
$$
R_{k_m}=N_{k_m}^\top M^{-1}{\mathcal P}_{{\mathcal A}_{k_m}}^\top({\mathcal {\mathcal P}_{{\mathcal A}_{k_m}}}M^{-1}{\mathcal P}_{{\mathcal A}_{k_m}}^\top)^{-1}
$$
and
$$
S_{k_m}=A_{k_m}+\frac{1}{\beta}N_{k_m}^\top M^{-1}N_{k_m}-\frac{1}{\beta}N_{k_m}^\top M^{-1} {\mathcal P}_{{\mathcal A}_{k_m}}^\top ({\mathcal P}_{{\mathcal A}_{k_m}} M^{-1}{\mathcal P}_{{\mathcal A}_{k_m}}^\top)^{-1}{\mathcal P}_{{\mathcal A}_{k_m}} M^{-1} N_{k_m}
$$
is the Schur complement of \eqref{RB_matrix}.

It is just seen that the Schur complement $S_{k_m}$ requires computing $({\mathcal P}_{{\mathcal A}_{k_m}} M^{-1}{\mathcal P}_{{\mathcal A}_{k_m}}^\top)^{-1}$. Note that the dimension of $M^{-1}$ is of order $O(h^{-d})$. Therefore, it is extremely expensive for fine mesh cases to compute $M^{-1}$ and hence the Schur complement $S_{k_m}$. To tackle this issue, we consider the lumped mass matrix (see \cite{LarsonB:2013}) to approximate the mass matrix $M$, which is a diagonal matrix with the row sums of the mass matrix $M$ on the diagonal. That is, we have
\begin{equation}\label{eq:lmW}
  W_{i,i}:=\sum_{j=1}^{N}M_{i,j}=\sum_{j=1}^{N}(\phi_h^j,\phi_h^i),~~i=1,2,\cdots,N,
\end{equation}
with $\phi_h^i$ being the finite element basis functions in $V_h$.
Then the inverse of the diagonal matrix $W$ is easy to compute.
Accordingly, the block matrix
$$
\hat{F}^{k_m}:=\begin{pmatrix}
A_{k_m} &N_{k_m}^\top&0\\
-N_{k_m} &\beta W&{\mathcal P}_{{\mathcal A}_{k_m}}^\top\\
0&{\mathcal P}_{{\mathcal A}_{k_m}}&0\end{pmatrix}$$
is an approximation of $F^{k_m}$ and its block diagonal decomposition is
\begin{equation}\label{bd_F_hat}
  \hat{F}^{k_m}=
    \underbrace{\begin{pmatrix}
      I &\frac{1}{\beta}N_{k_m}^\top W^{-1}&G_{k_m}\\
      0&I&0\\
      0&\frac{1}{\beta}{\mathcal P}_{{\mathcal A}_{k_m}}W^{-1}&I\end{pmatrix}}_{L_{k_m}}
    \underbrace{\begin{pmatrix}
      H_{k_m}&0&0\\
      0 &\beta W&0\\
      0&0&-\frac{1}{\beta}({\mathcal P}_{{\mathcal A}_{k_m}}W^{-1}{\mathcal P}_{{\mathcal A}_{k_m}}^\top)\end{pmatrix}}_{C_{k_m}}
    \underbrace{\begin{pmatrix}
      I &0&\\
      -\frac{1}{\beta} W^{-1}N_{k_m}&I&\frac{1}{\beta}W^{-1}{\mathcal {\mathcal P}_{{\mathcal A}_{k_m}}}\\
      -G_{k_m}^\top&0&I\end{pmatrix}}_{R_{k_m}},
\end{equation}
where
$ G_{k_m}=N_{k_m}^\top W^{-1}{\mathcal P}_{{\mathcal A}_{k_m}}^\top({\mathcal {\mathcal P}_{{\mathcal A}_{k_m}}}W^{-1}{\mathcal P}_{{\mathcal A}_{k_m}}^\top)^{-1}
$
and
$$
H_{k_m}=A_{k_m}+\frac{1}{\beta} N_{k_m}^\top W^{-1}N_{k_m}-\frac{1}{\beta} N_{k_m}^\top W^{-1} {\mathcal P}_{{\mathcal A}_{k_m}}^\top ({\mathcal P}_{{\mathcal A}_{k_m}} W^{-1}{\mathcal P}_{{\mathcal A}_{k_m}}^\top)^{-1}{\mathcal P}_{{\mathcal A}_{k_m}} W^{-1} N_{k_m}.
$$
Also, it is easy to verify that the matrix ${\mathcal P}_{{\mathcal A}_{k_m}}W^{-1}{\mathcal P}_{{\mathcal A}_{k_m}}^\top$ is diagonal. Thus, it is easy to compute  $({\mathcal P}_{{\mathcal A}_{k_m}}W^{-1}{\mathcal P}_{{\mathcal A}_{k_m}}^\top)^{-1}$ and it holds that
$$
W^{-1}- W^{-1} {\mathcal P}_{{\mathcal A}_{k_m}}^\top ({\mathcal P}_{{\mathcal A}_{k_m}} W^{-1}{\mathcal P}_{{\mathcal A}_{k_m}}^\top)^{-1}{\mathcal P}_{{\mathcal A}_{k_m}} W^{-1}= \Pi_{{\mathcal I}_{k_m}} W^{-1}\Pi_{{\mathcal I}_{k_m}}.
 $$
Therefore, computational cost for the explicit formulation of the matrix $G_{k_m}$ is negligible and the matrix $H_{k_m}$ can be simplified as
\begin{equation}\label{Simplified_Hm}
H_{k_m}=A_{k_m}+\frac{1}{\beta} N_{k_m}^\top \Pi_{{\mathcal I}_{k_m}} W^{-1}\Pi_{{\mathcal I}_{k_m}} N_{k_m}.
\end{equation}
Thus, computing $M^{-1}$ is not required for the Schur complement $H_{k_m}$ in \eqref{Simplified_Hm} and it becomes easy to compute the block factorization \eqref{bd_F_hat} of $\hat{F}^{k_m}$. These features suggest us to relax the Newtonian system \eqref{Newton_a0} to a linear system with the coefficient matrix $\hat{F}^{k_m}$.

Note that the Newtonian system \eqref{Newton_a0} can be rewritten as
\begin{equation}\label{Newton_temp}
\hat{F}^{k_m}\begin{pmatrix}
\bm{r}^{k_{m+1}}\\ \bm{q}^{k_{m+1}}-\bm{q}^{k_m}\\(\bm{\eta}^{k_{m+1}})_{{\mathcal A}_{k_m}}\end{pmatrix}=\begin{pmatrix}\bm{d}_1\\\bm{d}_2\\\bm{d}_3\end{pmatrix}+\begin{pmatrix}
0\\ \beta (W-M)(\bm{q}^{k_{m+1}}-\bm{q}^{k_m})\\0\end{pmatrix}.
\end{equation}
The equation (\ref{Newton_temp}) is implicit because $\bm{q}^{k_{m+1}}$ appears in both sides. We consider a semi-implicit scheme for the Newtonian system \eqref{Newton_temp} by replacing $\bm{q}^{k_{m+1}}$ in the right-hand side with the known last outer iterate $\bm{q}^{k}$, and obtain
\begin{equation}\label{Newton_a1}
\hat{F}^{k_m}\begin{pmatrix}
\bm{r}^{k_{m+1}}\\ \bm{q}^{k_{m+1}}-\bm{q}^{k_m}\\(\bm{\eta}^{k_{m+1}})_{{\mathcal A}_{k_m}}\end{pmatrix}=\begin{pmatrix}\bm{\hat{d}}_1\\\bm{\hat{d}}_2\\\bm{\hat{d}}_3\end{pmatrix} ,~\text{where}~\begin{pmatrix}\bm{\hat{d}}_1\\\bm{\hat{d}}_2\\\bm{\hat{d}}_3\end{pmatrix}:=\begin{pmatrix}\bm{d}_1\\\bm{d}_2\\\bm{d}_3\end{pmatrix}+\begin{pmatrix}
0\\ \beta (W-M)(\bm{q}^{k}-\bm{q}^{k_m})\\0\end{pmatrix}.
\end{equation}
That is, for numerical implementation purpose, we relax the Newtonian system \eqref{Newton_a0} as the much easier linear system \eqref{Newton_a1}. Indeed, it follows from \eqref{bd_F_hat} that
$$
\begin{pmatrix}
\bm{r}^{k_{m+1}}\\ \bm{q}^{k_{m+1}}-\bm{q}^{k_m}\\(\bm{\eta}^{k_{m+1}})_{{\mathcal A}_{k_m}}\end{pmatrix}=R_{k_m}^{-1}C_{k_m}^{-1}L_{k_m}^{-1}\begin{pmatrix}\bm{\hat{d}}_1\\\bm{\hat{d}}_2\\\bm{\hat{d}}_3\end{pmatrix}.
$$
Hence, the procedure of solving the linear system \eqref{Newton_a1} via its Schur complement reduction can be summarized in Subroutine \ref{Solve_Newton_a}.

\begin{algorithm}[htpb]\caption{Solver for the Newtonian system \eqref{Newton_a1}.}\label{Solve_Newton_a}
  \begin{algorithmic}[1]
    \Function {Newton-Solver}{$\bm{\hat{d}}_1, \bm{\hat{d}}_2, \bm{\hat{d}}_3$}
    \State  $(\bm{\hat{d}}_1, \bm{\hat{d}}_2, \bm{\hat{d}}_3)^\top
      \leftarrow L_{k_m}^{-1}(\bm{\hat{d}}_1, \bm{\hat{d}}_2, \bm{\hat{d}}_3)^\top$.
    \State $(\bm{r}^{k_{m+1}}, \bm{\hat{d}}_2, \bm{\hat{d}}_3)^\top
      \leftarrow C_{k_m}^{-1}(\bm{\hat{d}}_1, \bm{\hat{d}}_2, \bm{\hat{d}}_3)^\top.$
    \State  $(\bm{r}^{k_{m+1}},\bm{q}^{k_{m+1}}-\bm{q}^{k_m}, (\bm{\eta}^{k_{m+1}})_{{\mathcal{A}}_{k_m}})^\top
      \leftarrow R_{k_m}^{-1}(\bm{r}^{k_{m+1}}, \bm{\hat{d}}_2, \bm{\hat{d}}_3)^\top$.
    \State \Return{$(\bm{q}^{k_{m+1}}-\bm{q}^{k_m}, (\bm{\eta}^{k_{m+1}})_{{\mathcal{A}}_{k_m}})$}.
		\EndFunction
  \end{algorithmic}\label{InADMM_algorithm}
\end{algorithm}

In Subroutine \ref{Solve_Newton_a}, the inverses of $L_{k_m}$ and $R_{k_m}$ are easy to compute as their permutation matrices are triangular. The computation of $C_{k_m}^{-1}(\bm{\hat{d}}_1, \bm{\hat{d}}_2, \bm{\hat{d}}_3)^\top$ requires solving the linear system
\begin{equation}\label{eq_Hm}
H_{k_m} \bm{r}=\bm{\hat{d}}_1,
\end{equation}
as well as computing $W^{-1}\bm{\hat{d}}_2$ and $ ({\mathcal P}_{{\mathcal A}_{k_m}}W^{-1}{\mathcal P}_{{\mathcal A}_{k_m}}^\top)^{-1} \bm{\hat{d}}_3$. Since both $W$ and ${\mathcal P}_{{\mathcal A}_{k_m}}W^{-1}{\mathcal P}_{{\mathcal A}_{k_m}}^\top$ are diagonal matrices, the computational load of Subroutine \ref{Solve_Newton_a} is dominated by solving \eqref{eq_Hm}. Note that \eqref{eq_Hm} is easy because its coefficient matrix $H_{k_m}$ is positive definite and its dimension is the same as that of \eqref{eq:A_{k_m}}, which is much less than that of \eqref{Newton_a0}. Thus, via Subroutine \ref{Solve_Newton_a}, the indefinite linear system \eqref{Newton_a0} in higher dimension is significantly alleviated.

\subsection{Optimization insights}\label{opt_ins}
As just shown, the indefinite linear system \eqref{Newton_a0} is relaxed to \eqref{Newton_a1} by its Schur complement reduction. Then it is interesting to analyze the corresponding relaxation of the underlying optimization problem and discern its difference from the desired one \eqref{qh_dis}. Indeed, the linear system \eqref{Newton_a1} is equivalent to:
\begin{equation}\label{Newton_temp_2}
F^{k_m}\begin{pmatrix}
\bm{r}^{k_{m+1}}\\ \bm{q}^{k_{m+1}}-\bm{q}^{k_m}\\(\bm{\eta}^{k_{m+1}})_{{\mathcal A}_{k_m}}\end{pmatrix}+\begin{pmatrix}
0\\ \beta(W-M)(\bm{q}^{k_{m+1}}-\bm{q}^k)\\0\end{pmatrix}=\begin{pmatrix}\bm{d}_1\\\bm{d}_2\\\bm{d}_3\end{pmatrix},
\end{equation}
which differs from the Newtonian system \eqref{Newton_a0} in the extra term $\beta(W-M)(\bm{q}^{k_{m+1}}-\bm{q}^k)$. Then, following the steps reversely in subsection \ref{ssec:acn}, it is easy to see that \eqref{Newton_temp_2} is exactly the corresponding Newtonian system if the active-set Newton method is applied to the following optimization problem:
  \begin{equation}\label{qh_proximal}
    \begin{aligned}
      &\min_{\bm{q}\in \mathbb{R}^N}\left\{J(\bm{q})+\frac{\beta}{2}\| \bm{q}-\bm{p}^k+\frac{\bm{\lambda}^k}{\beta}\|_{M}^2+\frac{\beta}{2}\| \bm{q}-\bm{q}^k\|_{W-M}^2\right\},\\
      &~~~\text{s.t.}\quad \bm{a}_0 \le \bm{q}\le \bm{a}_1.
    \end{aligned}
  \end{equation}
Since each entry of the mass matrix $M$ is positive \cite{LarsonB:2013}, together with the definition of $W$ in \eqref{eq:lmW}, it is easy to verify that $W-M$ is positive semidefinite. Then, the problem \eqref{qh_proximal} is still convex and the solution of \eqref{Newton_a1} converges to the solution of the problem \eqref{qh_proximal}. Note that we slightly abuse the notation and define $\frac{\beta}{2}\| \bm{q}-\bm{q}^k\|_{W-M}^2:=(\bm{q}-\bm{q}^k)^\top(W-M)(\bm{q}-\bm{q}^k)$ in (\ref{qh_proximal}), despite that $W-M$ is positive semidefinite.
 Furthermore, because of the equivalence between the Euclidean space and the space $V_h$, \eqref{qh_proximal} can be rewritten as
  \begin{equation}\label{qh_pro}
    \min_{q_h\in K_h}\left\{J(q_h)+\frac{\beta}{2}\| q_h-p_h^k+\frac{\lambda_h^k}{\beta}\|_{L^2(\Omega)}^2+\frac{\beta}{2}\|q_h-q_h^k\|_T^2\right\}.
  \end{equation}
In (\ref{qh_pro}), the semi-norm
$$
\|q_h\|_T^2:=\sum_{\tau\in\mathcal{T}_h}Q_{\tau,h}(q_h^2)-(q_h,q_h) \;\hbox{with} \; Q_{\tau,h}(g)=\frac{|\tau|}{d+1}\sum_{j=1}^{d+1}g(x_j^\tau)
$$
and
$\{x_j^\tau\}_{j=1}^{d+1}$ are the vertices of the $d$-simplex $\tau\in\mathcal{T}_h$, and $\mathcal{T}_h$ is a regular partition of $\Omega$. Hence, our numerical technique for tackling the difficult $q_h$-subproblem via solving \eqref{Newton_a1} can be represented as replacing the problem \eqref{eq:qh} with \eqref{qh_pro}, in which the objective function is regularized by a semi-proximal regularization term. In other words, applying the active-set Newton method along with the Schur complement reduction can be explained as replacing the optimization problem (\ref{eq:qh}) by the proximally regularized one (\ref{qh_proximal}). Replacing the $q_h$-subproblem \eqref{eq:ADMM_a_origin_a} with \eqref{qh_pro} in the ADMM (\ref{eq:ADMM_a_origin}) hence results in the so-called proximal ADMM, which has been well studied in the optimization area. We refer to, e.g., \cite{HeB:2002, HeB:2012}, for convergence of various proximal versions of the ADMM.

\begin{remark}
For the implicit equation \eqref{Newton_temp}, we can alternatively consider replacing the unknown $\bm{q}^{k_{m+1}}$ in the right-hand of \eqref{Newton_temp} with the last inner iterate $\bm{q}^{k_m}$, instead of the last outer iterate $\bm{q}^{k}$. The resulting semi-implicit equation remains the coefficient matrix $\hat{F}^{k_m}$ and the right-hand side in \eqref{Newton_a0}. In our numerical experiments, we use the warm start technique, meaning the initial iterate is set as $\bm{q}^{k_0}:=\bm{q}^{k}$, and as to be shown in numerical results, usually each inner loop only requires executing the active-set Newton method by one iteration. Hence, using $\bm{q}^{k_m}$ or $\bm{q}^{k}$ makes very little difference numerically. On the other hand, an advantage of using $\bm{q}^{k}$ in \eqref{Newton_temp} is that the resulting scheme can be theoretically explained as a proximal version of the ADMM with well known theoretical results as studied in the optimization area.
\end{remark}

\subsection{Implementation of the active-set Newton method for the $q_h$-subproblem \eqref{Newton_a1}}\label{sec:NtsA_III}

Now, we present the active-set Newton method for solving the $q_h$-subproblem  \eqref{Newton_a1} in Subroutine \ref{asNewton}.

\begin{algorithm}[htpb]\caption{An active-set Newton method for the $q_h$-subproblem.}\label{asNewton}
  \begin{algorithmic}[1]
    \Function{asNewton}{$\bm{q}^k,\bm{p}^k,\bm{\lambda}^k$}
      \State Set initial values: $0\leftarrow m$, $(\bm{q}^{k_m},J'(\bm{q}^{k_m}))\leftarrow (\bm{q}^k, J'(\bm{q}^k))$; $\epsilon>0$; $\text{``Tol''}>0$.
      \While{$m\le\text{MaxIter}$}
        \State Compute the active and inactive indices: ${\mathcal A}_{k_m}^+,{\mathcal A}_{k_m}^-,{\mathcal A}_{k_m},{\mathcal I}_{k_m}$.
        \State $\begin{pmatrix}\bm{\hat{d}}_1\\\bm{\hat{d}}_2\\\bm{\hat{d}}_3 \end{pmatrix}\leftarrow \begin{pmatrix}0\\  -{J'}(\bm{q}^{k_m})-\beta W(\bm{q}^{k_m}-\bm{q}^k)-\beta M(\bm{q}^k-\bm{p}^k+\frac{\bm{\lambda}^k}{\beta})\\{\mathcal P}_{{\mathcal A}_{k_m}^+}\bm{a}_1+{\mathcal P}_{{\mathcal A}_{k_m}^-}\bm{a}_0-{\mathcal P}_{{\mathcal A}_{k_m}}\bm{q}^{k_m}\end{pmatrix}$.
        \State $(\bm{\eta}^{k_{m+1}})_{{\mathcal{I}}_{k_m}}\leftarrow 0$.
        \State $(\bm{q}^{k_{m+1}}-\bm{q}^{k_m}, (\bm{\eta}^{k_{m+1}})_{{\mathcal{A}}_{k_m}})\leftarrow$
        \Call{Newton-Solver}{$\bm{\hat{d}}_1, \bm{\hat{d}}_2, \bm{\hat{d}}_3$}, (See Subroutine \ref{Solve_Newton_a}).
        \State $\bm{q}^{k_{m+1}}=\max(\epsilon,\bm{q}^{k_{m+1}} )$.
        \State $m\leftarrow m+1$.
        \State $J'(\bm{q}^{k_m})\leftarrow$\Call{Gradient}{$\bm{q}^{k_m}$}, (See Subroutine \ref{Compute_J'a}).
        \If{$\|F(\bm{q}^{k_m},\bm{\eta}^{k_m})\|<$ Tol}
          \State \Return $\bm{q}^{k+1}\leftarrow \bm{q}^{k_m}$, $\textbf{break}.$
        \EndIf
      \EndWhile
    \EndFunction
  \end{algorithmic}
\end{algorithm}

For each step of Subroutine \ref{asNewton}, we need to solve two linear systems: \eqref{eq:A_{k_m}} in Subroutine \ref{Compute_J'a} and \eqref{eq_Hm} in Subroutine \ref{Solve_Newton_a}. Since both of these linear systems are positive definite, we can use the preconditioned conjugate gradient (PCG) method to solve them. As \eqref{eq:A_{k_m}} is a discretized formulation of the elliptic equation \eqref{eq:IDCEE}, a popular way to construct the preconditioner $\tilde{A}_{k_m}$ is using the multigrid (MG) method (see, e.g. \cite{BriggsW:2000}), which uses the MG V-cycles associated with its coefficient matrix $A_{k_m}$ to approximate $A_{k_m}^{-1}$. For \eqref{eq_Hm} in Subroutine \ref{Solve_Newton_a}, we still use the preconditioner $\tilde{A}_{k_m}$ of \eqref{eq:A_{k_m}} for solving \eqref{eq_Hm}, though the MG V-cycles associated with $H_{k_m}$ may be closer to $H_{k_m}^{-1}$. Note that we do not use the  MG V-cycles associated with $H_{k_m}$ because it is expensive to compute the explicit formulation of $H_{k_m}$ and MG V-cycles require more computation if the explicit formulation of $H_{k_m}$ is unknown.

\begin{remark}
The update of $\bm{q}^{k_{m+1}}$ in Line 8 of Subroutine \ref{asNewton} is to ensure that $\bm{q}$ is positive. The norm $\|F(\bm{q},\bm{\eta})\|$ is defined as
$$
\max ({\|\text{error}}_1\|_{W^{-1}}, \|{\text{error}}_2\|)
$$
with
$$
  \left \{
  \begin{aligned}
    {\text{error}}_1&:=-{J'}(\bm{q})-\beta W(\bm{q}-\bm{q}^k)-\beta M(\bm{q}^k-\bm{p}^k+\frac{\bm{\lambda}^k}{\beta})+\bm{\eta} , \\{\text{error}}_2&:=\bm{\eta}-\max \{0, \bm{\eta}+c(\bm{q}-\bm{a}_1)\}-\min\{0, \bm{\eta}+c(\bm{q}-\bm{a}_0)\}.
  \end{aligned}\right.
$$
\end{remark}

\begin{remark}\label{rem_convergence}
Since the convexity of the functional $J(\cdot)$ is kept, the active-set Newton method summarized as Subroutine \ref{asNewton} is guaranteed to be convergent to a solution of \eqref{qh_dis}. We refer to, e.g., \cite{HinzeP:2008}, for rigorous analysis. Also, as analyzed in  \cite{HinzeP:2008}, Subroutine \ref{asNewton} is superlinearly convergent.
\end{remark}

\begin{remark}
It is clear that the computational cost of Subroutine \ref{asNewton} is dominated by computing the gradient $J'(\bm{q}^{k_m})$ (Subroutine \ref{Compute_J'a}) and the Newton step (Subroutine \ref{Solve_Newton_a}).
To compute the gradient $J'(\bm{q}^{k_m})$, its main computation is computing the discrete matrix $A_{k_m}$ and then solving the discretized elliptic equation \eqref{eq:A_{k_m}}. For the Newton step, its main computation is calculating the discrete matrix $N_{k_m}$ and then solving the linear system \eqref{eq_Hm}. As just analyzed, the linear system \eqref{eq_Hm} is positive definite and its dimension is the same as that of \eqref{eq:A_{k_m}}. Hence, the Newton step does not require too much additional computation, compared with the computation of the gradient $J'(\bm{q}^{k_m})$. Recall that implementing a first-order algorithm usually requires discerning an appropriate step size (e.g., via line-search techniques) for the sake of ensuring the convergence, hence multiple objective function values are usually required. As mentioned, computing these functional values is equally expensive as that of computing the gradient. Therefore, it is encouraging to consider the active-set Newton method in Subroutine \ref{asNewton} whose computation is not much more than that of implementing a first-order algorithm, yet its convergence is guaranteed to be superlinear.
\end{remark}

\section{Deep CNN for the $p_h$-subproblem}\label{sec:CNN}

In this section, we discuss how to solve the $p_h$-subproblem (\ref{eq:ADMM_a_origin_b}). This subproblem can be specified as
\begin{equation}\label{eq:ph_TV}
  p_h^{k+1}=\underset{p_h\in V_h}{\arg\min}\;\Big\{ \alpha\| \nabla p_h\|_{L^1(\Omega)}+\frac{\beta}{2}\| q_h^{k+1}-p_h+\frac{\lambda_h^k}{\beta}\|_{L^2(\Omega)}^2\Big\}.
\end{equation}
Note that the original TV term is kept and hence the objective functional in (\ref{eq:ph_TV}) is nonsmooth. Obviously, \eqref{eq:ph_TV} has no closed-form solution and it should be solved iteratively by a certain algorithm. Also, the dimension of $p_h$ is the same as that of $q_h$, and it may be high for a higher-dimensional space and fine mesh discretization. For instance, it is of order $10^6$ if the mesh size $h=1/1024$ for the unit square $\Omega \subset \mathbb{R}^2$. Hence, it is also necessary to consider how to solve the $p_h$-subproblem (\ref{eq:ADMM_a_origin_b}) efficiently for implementing the ADMM (\ref{eq:ADMM_a_origin}). We reiterate that it is always more preferable to choose some model-tailored algorithms in accordance with the structure of the problem under consideration. For the $p_h$-subproblem (\ref{eq:ph_TV}), certainly it can be treated as a generic optimization problem and then some generic-purpose or less structure-exploiting algorithms can be applied. But it turns out that the deep convolutional neural network (CNN), which has been significantly enhanced in recent literatures (e.g., \cite{KrizhevskyA2012,SimonyanK2015}), is a much better choice for the $p_h$-subproblem (\ref{eq:ph_TV}). Below is the detail.

Let $R(x):=\| \nabla p_h\|_{L^1(\Omega)}$ and $\theta>0$ be constant. The proximal operator of $R(x)$ is given by
\begin{equation}\label{denoising_operator}
  \text{Prox}_{\theta R(x)}(z)=\arg \min_x \Big\{\theta  R(x)+\frac{1}{2}\|x-z\|_{L^2(\Omega)}^2\Big\}.
\end{equation}
Then, the solution of the $p_h$-subproblem \eqref{eq:ph_TV} can be presented by
$$
  p_h^{k+1}=\text{Prox}_{\frac{\alpha}{\beta}R(x)}(q_h^{k+1}+\frac{\lambda_h^k}{\beta}).
$$
Following the standard Rudin-Osher–Fatemi model in \cite{RudinOF:1992}, the operator $\text{Prox}_{\frac{\alpha}{\beta} R(x)}$ can be interpreted as the denoising operator for the standard image denoising model. In the last few years, the literature of algorithms for various image denoising models has been phenomenally upgraded by contemporary deep neural networks, see, e.g.,\cite{xie2012image,ZhangZCMZ:2017}. An advantage of applying a deep neural network to denoising models is that it avoids iterations in its testing phase,
and hence computation can be largely saved. We are thus inspired to consider some pre-trained deep neural network, rather than some iterative scheme, for the $p_h$-subproblem (\ref{eq:ph_TV}).

To see why the deep CNN is chosen for the case where $\Omega\subset \mathbb{R}^2$ is a rectangular domain and it is triangulated into the uniform mesh, there exists a one-to-one mapping between $p_h\in V_h$ and an $m\times n$ raster image ($m\times n=N$) where the gray value at pixel $(i,j)$ of the image $I$ corresponds to the value of the function $p_h$ at node $(i,j)$. Thus, there is a mapping between a discrete two-dimensional function and a gray-scale raster image. Then, the pre-trained deep CNN which has been widely used for various image denoising problems can be applied.
Let $\mathcal{M}$ denote the mapping from $\bm{p}$ to a raster image, and $\mathcal{C}_{\sigma}$ the pre-trained deep CNN with $\sigma$ the variance of the noise used for training CNN. Solving the $p_h$-subproblem (\ref{eq:ph_TV}) by a pre-trained deep CNN can be summarized in Subroutine \ref{Algo_7}.
\begin{algorithm}[htpb]\caption{A deep CNN based method for the $p_h$-subproblem.}\label{Algo_7}
  \begin{algorithmic}[1]
    \Function{Dnoiser}{$\bm{q}^{k+1},\bm{\lambda}^k$}
      \State $I^{k+1}_{input}:=\mathcal{M}(\bm{q}^{k+1}+\frac{\bm{\lambda}^k}{\beta})$.
      \State $I^{k+1}_{output}:=\mathcal{C}_{\sigma }(I_{input}^{k+1})$.
      \State $\bm{p}^{k+1}=\mathcal{M}^{-1}(I^{k+1}_{output})$.
    \EndFunction
  \end{algorithmic}
\end{algorithm}

\begin{remark}
Our primary interest is the case where $\Omega\subset \mathbb{R}^2$ is a rectangular domain and it is partitioned by the uniform triangulation mesh. For other cases such as $\Omega$ is not rectangular, the mesh is not uniform, or the deep CNN is not trained based on raster images, the mapping $\mathcal{M}$
should be redefined. For the case where $\Omega \subset \mathbb{R}^3$, one may employ a
deep 3D CNN (see, e.g., \cite{Ji2012}). These much more complicated situations should be discussed case by case with significantly more techniques, and they are beyond the scope of this paper.
\end{remark}

\section{The ADMM-Newton-CNN numerical approach}\label{sec:ADMM-Newton-CNN}

\begin{figure}[H]
  \centering
  \begin{tikzpicture}[node distance=1cm]
    \node (Input) [rectangle_rounded_green, yshift=0.7cm] {Input $\nabla u_{\delta}$, $f$, $a_0$, $a_1$; compute: $M$, $W$; set initials $\bm{q}^0$, $\bm{u}^0$, $\bm{\eta}^0$, $\bm{p}^0$, $\bm{\lambda}^0$; set $k=0$, $m=0$.};
    \node (compute_A) [rec, below of=Input, yshift=-0.5cm] {Compute $A_{k_m}$};
    \node (compute_u) [rec, right of=compute_A, xshift=2cm] {Compute $\bm{u}^{k_m}$};
    \node (compute_J') [rec, below of=compute_u, xshift=3cm] {Compute $J'$};
    \node ( Subroutine 1) [rec_lable, right of=compute_u, xshift=2cm ] {{\color{red}Subroutine 1}};
    \node ( Subroutine 3) [rec_lable, right of=compute_u, xshift=5cm] {{\color{red}Subroutine 3}};
    \node (Check_inner) [ellipse_red, below of=compute_J', yshift=-0.6cm] {Whether Newton method converges?};
    \node (compute_as) [rec, left of=Check_inner, xshift=-8cm] {Compute the active sets and  $N_{k_m}$};
    \node (compute_Newton) [rec_long, left of=compute_J', xshift=-8cm, yshift=-0.2cm] {{\small Newton step: call {\color{red}Subroutine 2} to get $(\bm{q}^{k_{m+1}}, (\bm{\eta}^{k_{m+1}})_{{\mathcal{A}}_{k_m}})$; set $ (\bm{\eta}^{k_{m+1}})_{{\mathcal{I}}_{k_m}}=0$.}};
    \node (Update_{k_m}) [trapezium_blue , left of=compute_A, xshift=-2cm] {$m\leftarrow m+1$ };
    \node (NO_inner) [rec_lable, left of=Check_inner, xshift=-4.5cm, yshift=0.2cm ] {No};
    \node (output_q) [trapezium_blue, below of=Check_inner, yshift=-1cm] {$\bm{q}^{k+1}\leftarrow\bm{q}^{k_m}$};
    \node (Yes_inner) [rec_lable, below of=Check_inner, xshift=-0.5cm] {Yes};
    \node (Check_outer) [ellipse_red, below of=output_q, yshift=-1cm] { Whether a satisfactory $\bm{q}$ is obtained?};
    \node (Update_k) [trapezium_blue,below of=compute_as, yshift=-1cm] {$k\leftarrow k+1$ };
    \node (reset_{k_m}) [trapezium_blue,below of=compute_as] {$m\leftarrow 0$ };
    \node (Update_lambda) [rec , below of=Update_k] {$\bm{\lambda}^{k+1}=\bm{\lambda}^k+\beta(\bm{q}^{k+1}-\bm{p}^{k+1})$};
    \node (Update_p) [rec , below of=Update_lambda] {Call {\color{red}Subroutine 4} to get $\bm{p}^{k+1}$};

    \node (output_final) [rectangle_out,below of=Check_outer, yshift=-0.5cm] {Output $\bm{q}^{k+1}$};
    \node (Yes_outer) [rec_lable, below of=Check_outer, xshift=-0.5cm, yshift=0.2cm ] {Yes};
    \node (NO_inner) [rec_lable, right of=Check_outer, xshift=-6.5cm, yshift=0.2cm ] {No};
    \draw [arrow](Input) -- (compute_A) ;
    \draw [arrow](compute_A) --(compute_u);
    \draw [arrow](compute_u) |-(compute_J');
    \draw [thick,draw=blue](-1.5,-0.3) -|(7.5,-2.5);
    \draw [thick,draw=blue](-1.5,-0.3) |-(2,-1.3);
    \draw [thick,draw=blue] (2,-1.3) |-(7.5,-2.5);
    \draw [arrow](Check_inner)--(compute_as) ;
    \draw [arrow](compute_as)--(compute_Newton) ;
    \draw [arrow](compute_Newton)--(Update_{k_m});
    \draw [arrow](Update_{k_m})--(compute_A);
    \draw [arrow](compute_J') --(Check_inner);
    \draw [arrow] (Check_inner)-- (output_q);
    \draw [thick,draw=blue](-6.8,-0.1) -|(10.5,-4.85);
    \draw [thick,draw=blue](-6.8,-0.11) |-(10.5,-4.85);
    \draw [arrow] (output_q)-- (Check_outer);
    \draw [arrow]  (Check_outer)--(Update_p);
    \draw [arrow]  (Update_p)--(Update_lambda);
    \draw [arrow]   (Update_lambda)--  (Update_k);
    \draw [arrow] (Update_k)-- (reset_{k_m});
    \draw [arrow] (reset_{k_m})--(compute_as);
    \draw [arrow]  (Check_outer)--(output_final);
  \end{tikzpicture}
  \caption{Flowchart of implementation of the proposed ADMM-Newton-CNN approach.}\label{Flowchart1}
\end{figure}
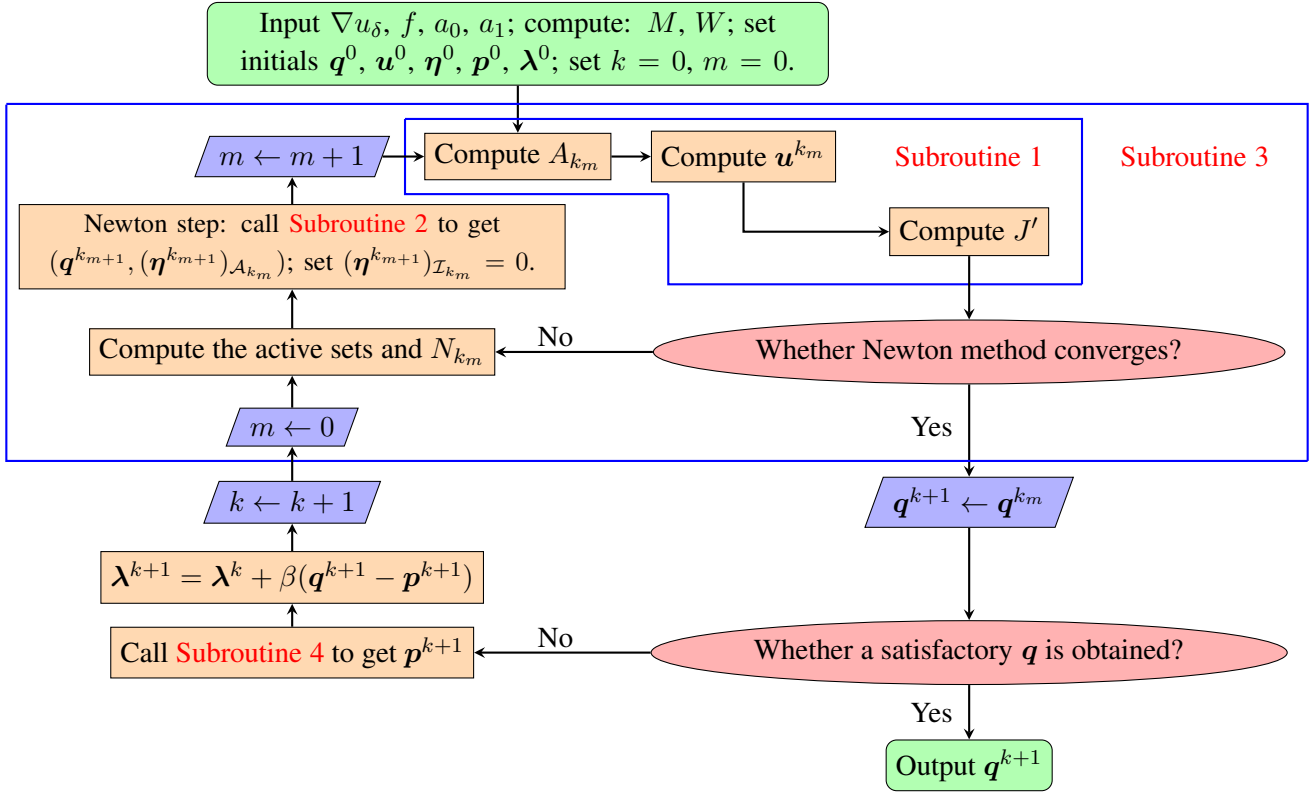

With the discussions in Sections \ref{sec:ADMM} and \ref{sec:CNN}, we are ready to present the complete version of the ADMM-Newton-CNN numerical approach to the TV model \eqref{eq:ROTV}. We show the flowchart of its implementation in Figure \ref{Flowchart1}.

\section{Numerical results}\label{sec:num}
In this section, we show efficiency of the proposed ADMM-Newton-CNN numerical approach by preliminary numerical results. All codes were written in MATLAB R2020b and numerical experiments were conducted on a desktop with Windows 10, Intel(R) Core(TM) i9-9900KF CPU (3.60 GHz), and 128 GB RAM. We notice that there are some limited numerical studies in the literature \cite{ChenZ:1999, Zou:1998}, which are focused on the smoothing $TV_\epsilon$ model \eqref{eq:dBMOP_2} with $d=1$ and coarse mesh discretization. But the proposed ADMM-Newton-CNN approach is for the original TV model \eqref{eq:ROTV} with the focus on the higher-dimensional space of $d=2$ as well as fine mesh discretization. Hence, it seems difficult to make any numerical comparison with the mentioned existing works, because of the essentially different natures in both modeling and philosophy of algorithmic design.

\subsection{Experiment setups}
We fix $\Omega$ as $(0,1)\times (0.1)$ and $f=10$ in $\Omega$.
The domain $\Omega$ is partitioned by the uniform triangulation mesh in the iFEM package \cite{iFEM}. The lower and upper bounds $a_0$ and $a_1$ in the constrained set $K$ in \eqref{eq:K} are taken as $0.1$ and $5.0$, respectively.
We follow \cite{ChenZ:1999,KeungZ:2001} and construct examples for the test in the following way.
\begin{enumerate}
  \item [(1.)] Choose a discontinuous diffusion coefficient $q(x)\in L^1(\Omega)$.
  \item [(2.)] Compute the finite element solution $u_h$ of \eqref{eq:IDCEE}.
  \item [(3.)] Take the noisy observation data as $   \nabla u_{\delta}(x)=\nabla u_h+\delta \*\|\nabla u_h\|_h \*\text{rand}(x)$,   where $\text{rand}(x)$ is a uniformly distributed random vector-valued function in $[-1,1]$ with $\delta>0$ the noise level.
\end{enumerate}
Recall that, for the $q_h$-subproblem \eqref{eq:ADMM_a_origin_a}, we use the warm start technique for the PCG executions, and the implementation of MG V-circles is based on the iFEM package developed in \cite{iFEM} with Jacobi splitting.
Moreover, for the $p_h$-subproblem \eqref{eq:ADMM_a_origin_b}, we use the pre-trained deep CNNs in \cite{ZhangZCMZ:2017} (\url{https://github.com/cszn/DnCNN}) and the mapping $\mathcal{M}$ in Subroutine \ref{Algo_7} is specified as $\mathcal{M}(\bm{p})=(\bm{p}-\bm{a}_0)/(a_1-a_0)\times256$, where $\bm{p}$ is the coefficient of $p_h$ and $I$ is a raster image. The initial guess of the Lagrange multiplier $\bm{\lambda}^0$ is always set to be $\bm{0}$; the initial guess of $\bm{q}^0$, $\bm{u}^0$ and $\bm{p}^0$ are set to be $\bm{1}$, $\bm{1}$ and $\bm{0}$, respectively.  For the stopping criterion to solve the linear systems  \eqref{eq:A_{k_m}}  and \eqref{eq_Hm}, relative errors are controlled with the tolerances of $10^{-10}$ and $10^{-5}$, respectively.
In addition, the value of $\text{``Tol''}$ in Subroutine \ref{asNewton} is $10^{-3}$.

\subsection{Experimental results}
\begin{example}\label{exm1}
  We take the discontinuous coefficient $q(x,y)$ in $(0,1)\times(0,1)$ as
  $$
    q(x,y)=
    \left\{
    \begin{aligned}
      &1, ~~ y\in [0,0.5],\\
      &2, ~~ y\in (0.5,1],\\
    \end{aligned}
    \right.
  $$
  whose discontinuous points form a straight line.
\end{example}

For the penalty parameter $\beta$ and denoising parameter $\sigma$, generally they should be tuned according to the noise level $\delta$. According to the Morozov's discrepancy principle \cite{EnglHN:1996},  the value of $\alpha$ in \eqref{eq:ROTV} is positively correlated with the noise level $\delta$, and the parameters $\sigma$ in $\mathcal{C}_{\sigma }$ and $\theta=\frac{\alpha}{\beta}$ in \eqref{denoising_operator} play the same role of controlling the rate of denoising. Hence, $\sigma$ should be proportional to $\theta=\frac{\alpha}{\beta}$ and $\beta \sigma$ should be positively correlated with the noise level $\delta$. In our numerical experiments, we tune the parameters $\sigma$ and $\beta$ such that $\beta \sigma$ is proportional to the noise level of the observation, i.e., $\beta \sigma \sim \delta$. In Table \ref{parameter_ex1}, we list the tuned values of $\beta$ and $\sigma$ for the cases where the noise levels are $\delta=0.01$, $0.05$ and $0.1$, respectively. These parameters are kept as constants for different finite element meshes.

\begin{table}[htb]
  \centering
  \caption{Parameters $\delta$, $\beta$ and $\sigma$ for Example \ref{exm1}.}\label{parameter_ex1}
    \vspace{3pt}
    \setlength{\tabcolsep}{18pt}\small
    \begin{tabular}{|c|c|c|c|}
	  \hline
      $\delta$ &$\beta$&$\sigma$&$(\beta\sigma)\slash\delta$\\
	  \hline
	  $0.01$ &0.1&9 &90\\
	  \hline
	  $0.05$ &0.5&9&90\\
	  \hline
	    $0.1$ &0.6&15&90\\
	  \hline
  \end{tabular}
\end{table}

\begin{table}[htb]
  \centering
  \caption{Numerical results for Example \ref{exm1} after the first 50 iterations.}\label{mesh_a_ex}
  \vspace{3pt}\setlength{\tabcolsep}{5pt}
  \begin{tabular}{|c|c|c|c|c|c|c|}
    \hline
    $\delta$&$h$&Total Newton No.&Total PCG No. for \eqref{eq:A_{k_m}}/\eqref{eq_Hm}&CPU Time (s) & $\|q_h^{50}-q\|_{L^2(\Omega)}/\|q\|_{L^2(\Omega)}$\\
   \hline
    &1/64&   57&   661 / 2238&2.425& 0.0068\\
       \cline{2-6}
    &1/128& 58&705 / 2722&9.327&0.0055\\
     \cline{2-6}
   0.01 &1/256&60&764 /  3239&45.847& 0.0053\\
 \cline{2-6}
    &1/512&62&814 / 3635&199.392& 0.0044\\
     \cline{2-6}
    &1/1024&63&855 / 3825&881.665& 0.0052\\
    \hline
    \hline
        &1/64&   55& 636 / 1211&2.093& 0.0106\\
       \cline{2-6}
    &1/128&55&660 / 1388&7.550&0.0122\\
     \cline{2-6}
   0.05 &1/256&55&679 / 1533&35.839& 0.0139\\
 \cline{2-6}
    &1/512&55&705 /1659&148.891& 0.0273\\
     \cline{2-6}
    &1/1024&55&724 / 1705& 649.436& 0.0486\\
    \hline
    \hline
        &1/64&  55&  638 / 1075&2.063& 0.0405\\
       \cline{2-6}
    &1/128& 55&666 / 1232&7.363&0.0374\\
     \cline{2-6}
   0.1 &1/256&55&683 / 1348&34.702& 0.0369\\
 \cline{2-6}
    &1/512&55&716 / 1448&144.838& 0.0478\\
     \cline{2-6}
    &1/1024&55&741 / 1500&631.207& 0.0706\\
    \hline
  \end{tabular}
\end{table}
\begin{table}[htb]
  \centering
  \caption{Computing time of various subtasks for Example \ref{exm1} with $h=1/256$  after the first 50 iterations.}\label{performance_ex}
  \vspace{3pt}
  \setlength{\tabcolsep}{3pt}\small
  \begin{tabular}{|c|c|c|c|c|c|c|c|c|c|c|c|}
	\cline{1-4} \cline{5-8} \cline{9-12}
    &\multicolumn{3}{c|}{$\delta$=0.01}&&\multicolumn{3}{|c|}{$\delta$=0.05}&&\multicolumn{3}{|c|}{$\delta$=0.1}\\
\cline{1-4} \cline{6-8} \cline{10-12}
	Subtasks&No.&Total time (s)& $\%$Time&&No.&Total time (s)& $\%$Time&&No.&Total time (s)& $\%$Time\\
	\cline{1-4} \cline{6-8} \cline{10-12}
	 Linear system \eqref{eq_Hm} & 60&16.901&$36.9\%$&        & 55&8.300&$23.2\%$&     &55&7.178&$20.7\%$\\
\cline{1-4} \cline{6-8} \cline{10-12}
      Implementation of CNN &50&9.256&$20.8\%$&        &50&9.810&$27.4\%$&        &50&9.556&$27.5\%$\\
\cline{1-4} \cline{6-8} \cline{10-12}
 $N_{k_m}$&60&4.578 &$10.0\%$&      &55&4.099 &$11.4\%$&       &55&4.179 &$12.0\%$\\
\cline{1-4} \cline{6-8} \cline{10-12}
$A_{k_m}$&60&3.789&$8.3\%$&    &55&3.440&$9.6\%$&        &55&3.478&$10.0\%$\\
\cline{1-4} \cline{6-8} \cline{10-12}
Linear system \eqref{eq:A_{k_m}} &60&3.196&$7.0\%$&     &55&2.842&$7.9\%$&      &55&2.897&$8.3\%$\\
\cline{1-4} \cline{6-8} \cline{10-12}
Others &&7.857&$17.0\%$&                       &&7.347&$20.5\%$&    &&7.414&$21.5\%$\\
\cline{1-4} \cline{6-8} \cline{10-12}
Total&&45.847&$100\%$&                   &&35.839&$100\%$&       &&34.702&$100\%$\\
	\hline
  \end{tabular}
\end{table}

As discussed in subsection \ref{opt_ins} and Remark \ref{rem_convergence}, the proposed ADMM-Newton-CNN approach is guaranteed to be convergent and our main interest is to show how numerically efficient this scheme could be. We have observed that the iterative sequence tends to be convergent after about $30$ iterations. Hence, we record the numerical performance in Table \ref{mesh_a_ex} for the first 50 iterations. For succinctness, only several choices of the noise levels and the meshes are listed. It is encouraging to see that total numbers of Newton steps, and PCG numbers for solving the linear systems \eqref{eq:A_{k_m}} and \eqref{eq_Hm}, as well as the relative error to the true solution $\|q_h^k-q\|_{L^2(\Omega)}/\|q\|_{L^2(\Omega)}$, are all very robust to the mesh.
Since the dimension of the resulting subproblems is increased when the mesh is refined, this feature is particularly favorable for fine mesh discretization.

To take a closer look into computing time of individual subtasks, we focus on the case of $h=1/256$ and report the respective computing times of various subtasks of the first 50 iterations in Table \ref{performance_ex}. According to this table, we see that computing time for the linear system \eqref{eq_Hm} accounts for about $20$-$40 \%$ of the entire time. Especially, for the cases where $\delta=0.01$ and $0.05$, computing time for \eqref{eq_Hm} is less than that of the CNN implementation. This fact well explains that the preconditioner $\tilde{A}_{k_m}$ is a good choice for the linear system \eqref{eq_Hm}. The Newton step is hence computationally cheap because the linear system \eqref{eq_Hm} can be well solved with the preconditioner $\tilde{A}_{k_m}$.
Recall that the computation of both ${J'}(\bm{q}^{k_m})$ and the objective function value mainly consists of computing $A_{k_m}$ and solving the linear system \eqref{eq:A_{k_m}}. Also, the Newton step needs to compute $N_{k_m}$ and solve the linear system \eqref{eq_Hm}. Based on Table \ref{performance_ex}, it is easy to estimate that the computation time of the Newton step is only about three times of that of computing ${J'}(\bm{q}^{k_m})$. Hence, choosing the active-set Newton method in Subroutine \ref{asNewton} for \eqref{qh_dis}, instead of some first-order algorithm which generally requires computing the objective function values repeatedly to find an appropriate step size, is verified.

\begin{figure}[h]
  \center
  \includegraphics[width=0.9\textwidth]{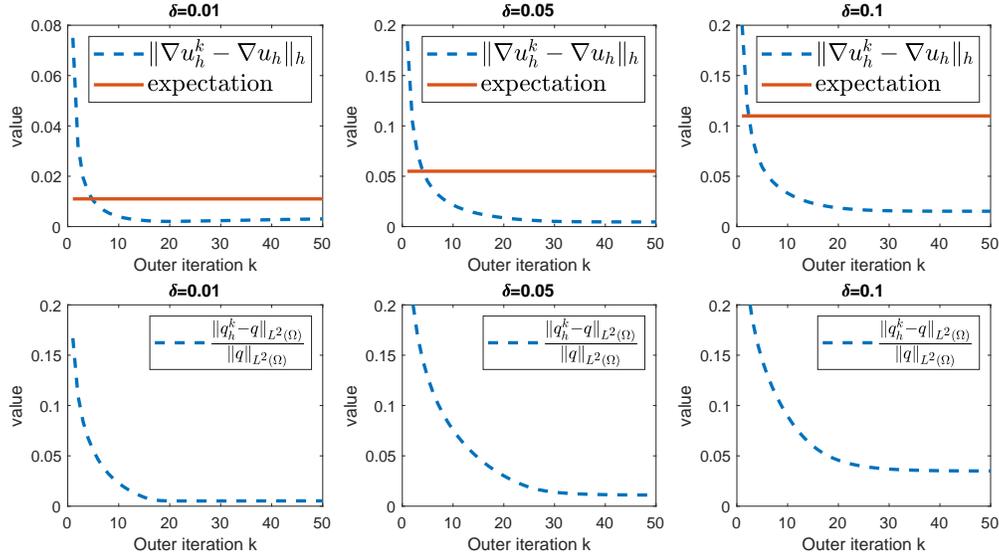}
  \caption{Plots of $\| \nabla u_h^k-\nabla u_h\|_h$ and $\|q_h^k-q\|_{L^2(\Omega)}/\|q\|_{L^2(\Omega)}$ for Example \ref{exm1} with $h=1/256$.}\label{iterative_result_ex1}
\end{figure}

\begin{figure}[h]
  \center
  \includegraphics[width=\textwidth]{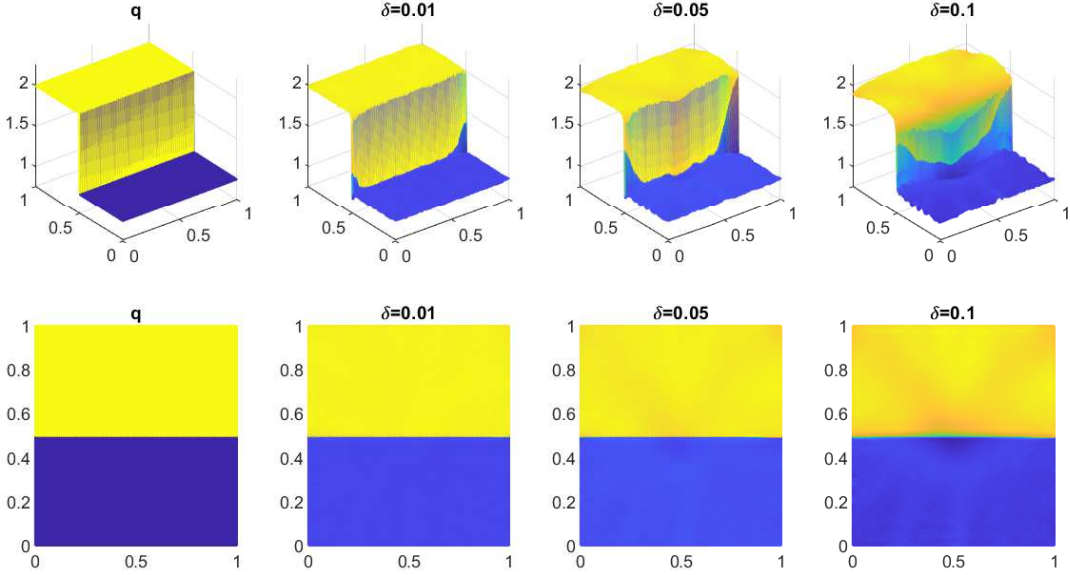}
  \caption{Numerical solutions $q_h^k$ at the 30-th iteration for Example \ref{exm1} with $h=1/256$. Column 1: the true coefficient $q$; Column 2: $\delta=0.01$ and the relative error $\|q_h^k-q\|_{L^2(\Omega)}/\|q\|_{L^2(\Omega)}=0.0051$; Column 3: $\delta=0.05$ and the relative error $\|q_h^k-q\|_{L^2(\Omega)}/\|q\|_{L^2(\Omega)}=0.0138$;  Column 4: $\delta=0.1$ and the relative error $\|q_h^k-q\|_{L^2(\Omega)}/\|q\|_{L^2(\Omega)}=0.0368$; Bottom: the projection of $q$ or $q_h^k$ onto the domain $\Omega$.}\label{iterative_result_ex1a_qh}
\end{figure}

In Figure \ref{iterative_result_ex1}, we plot the curves of $\| \nabla u_h^k-\nabla u_h\|_h$ and $\|q_h^k-q\|_{L^2(\Omega)}/\|q\|_{L^2(\Omega)}$ for Example \ref{exm1} with $h=1/256$. Also, differences between the ground-truth solution $q$ and the numerical solutions $q_h^k$ at the 30-th iteration are plotted in Figure \ref{iterative_result_ex1a_qh}. These curves further display the efficiency of the proposed ADMM-Newton-CNN approach for Example \ref{exm1}. In this figure, ``expectation'' means $E[\| \nabla u_h-\nabla u_h\|_h]$.

\begin{example}\label{exm2}
 The discontinuous coefficient $q(x,y)$ is taken as
  $$
    \begin{aligned}
      q(x,y)&=1+0.5*I_{\Omega_1}+I_{\Omega_2},\\ \Omega_1&=\{(x,y)|(x-0.5)^2+(y-0.5)^2\le1/8\},\\ \Omega_2&=\{(x,y)|1/3 \le x\le 2/3, 1/3 \le y\le 2/3\},
    \end{aligned}
  $$
  where $I_{\Omega_k}$ denotes the characteristic function over $\Omega_k,~k=1,2$. Its discontinuous points form a circle and a square.  This example has right-angled and curved discontinuous points, and it is more complicated.
\end{example}

Values of the parameters $\beta$ and $\sigma$ for various noise levels $\delta$ are listed in Table \ref{parameter_ex2}. Again, values of $\sigma$ are set such that $\beta \sigma$ is proportional to the noise level of the observation, i.e., $\beta \sigma \sim \delta$, and these parameters are kept as constants for different finite element meshes.

\begin{figure}[h]
  \center
  \includegraphics[width=0.75\textwidth]{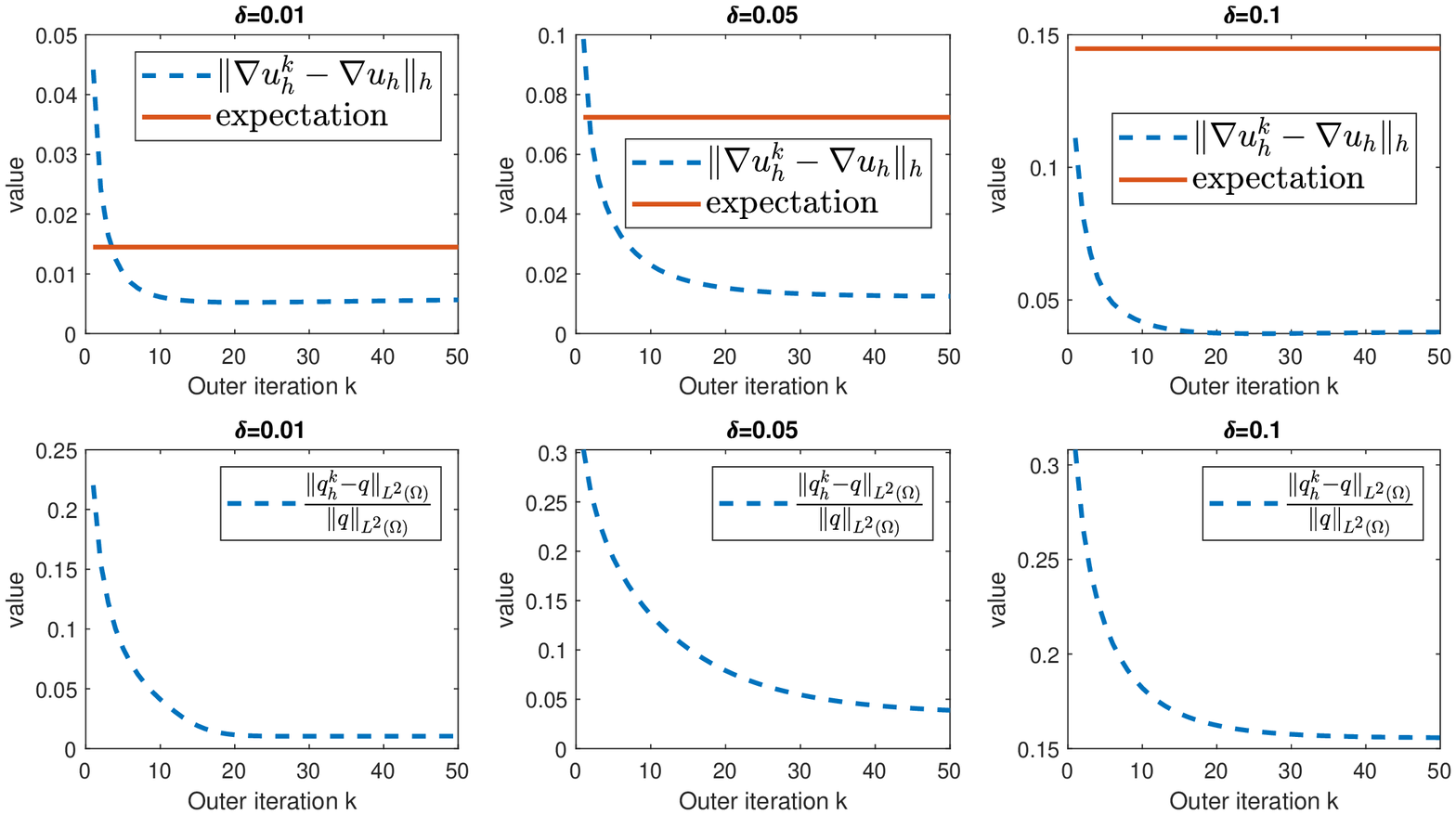}
  \caption{Plots of $\| \nabla u_h^k-\nabla u_h\|_h$ and $\|q_h^k-q\|_{L^2(\Omega)}/\|q\|_{L^2(\Omega)}$ for Example \ref{exm2} with $h=1/256$}\label{iterative_result_ex2}
\end{figure}

\begin{figure}[h]
  \center
  \includegraphics[width=\textwidth]{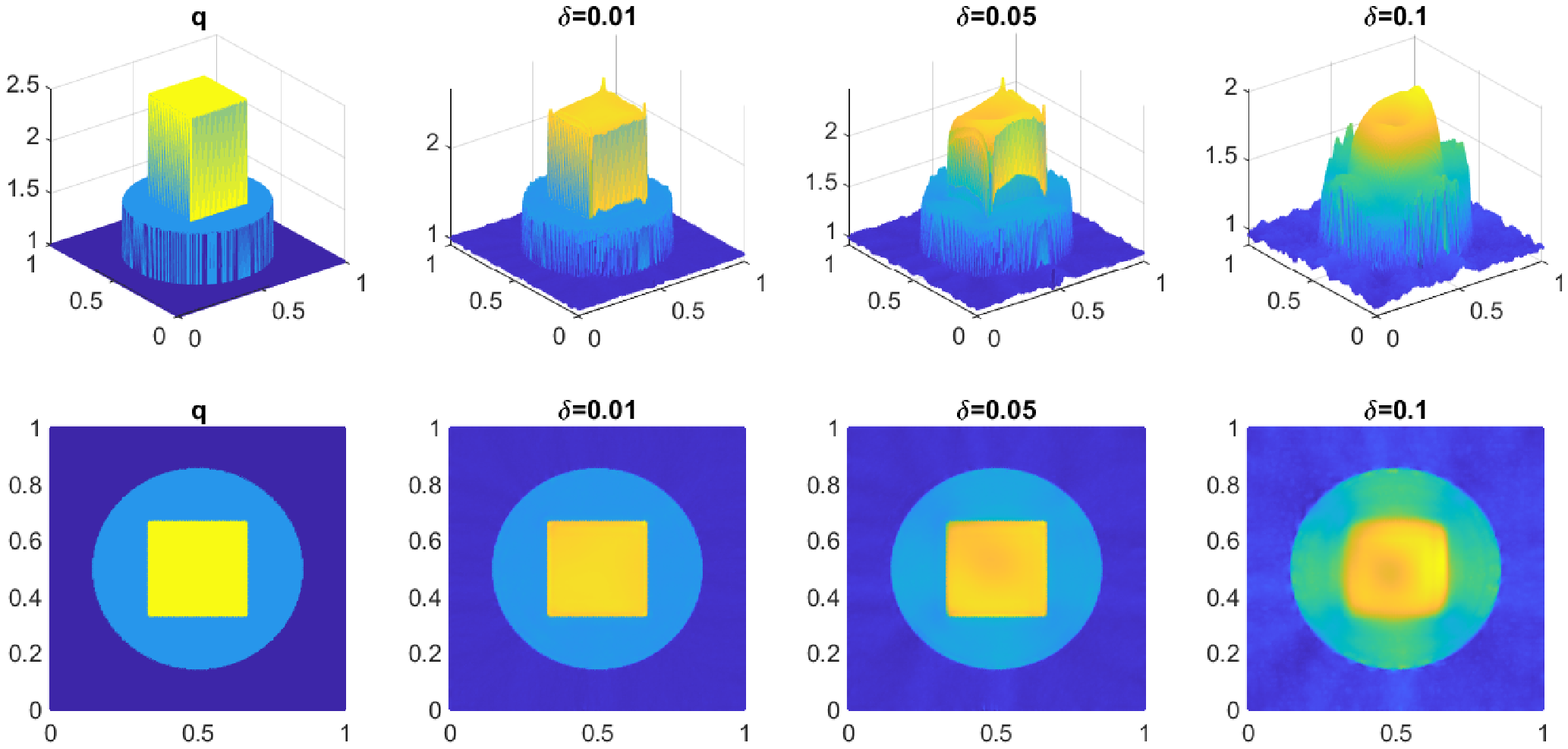}
  \caption{Numerical solutions $q_h^k$ at the 30-th iteration for Example \ref{exm2} with $h=1/256$. Column 1: the true coefficient $q$; Column 2: $\delta=0.01$ and the relative error $\|q_h^k-q\|_{L^2(\Omega)}/\|q\|_{L^2(\Omega)}=0.0104$; Column 3: $\delta=0.05$ and the relative error $\|q_h^k-q\|_{L^2(\Omega)}/\|q\|_{L^2(\Omega)}=0.0544$;  Column 4: $\delta=0.1$ and the relative error $\|q_h^k-q\|_{L^2(\Omega)}/\|q\|_{L^2(\Omega)}=0.1575$; Bottom: the projection of $q$ or $q_h^k$ onto the domain $\Omega$.}\label{iterative_result_ex2a_qh}
\end{figure}

Numerical results are reported in Table \ref{mesh_a_ex2}, for the first 50 iterations when the proposed ADMM-Newton-CNN approach is applied to Example \ref{exm2}. We list computing time of various subtasks individually in Table \ref{performance_ex2}, for the  first 50 iterations and $h=1/256$. Moreover, in Figure \ref{iterative_result_ex2}, $\|\nabla u_h^k-\nabla u_h\|_h$ and $\|q_h^k-q\|_{L^2(\Omega)}/\|q\|_{L^2(\Omega)}$ are plotted for the first 50 iterations when $h=1/256$. In this figure, ``expectation'' means $E[\| \nabla u_h^k-\nabla u_h\|_h]$. In Figure \ref{iterative_result_ex2a_qh}, difference between the ground-truth solution $q$ and  $q_h^k$ at the 30-th iteration are plotted for the case where $h=1/256$. Similar conclusions as those for Example \ref{exm1} can be drawn, and efficiency of the proposed ADMM-Newton-CNN approach is further verified for Example \ref{exm2}.

\section{Conclusions}\label{sec:Conc}
We focus on a well-known model with the total variational (TV) regularization for identifying the diffusion coefficient in an elliptic equation with observation data of the gradient of the solution. We consider the original TV-regularized model without any relaxation so that the favorable nonsmoothness and convexity properties can be both kept. We propose to solve this model by the alternating direction method of multipliers (ADMM), and show that the resulting subproblems can be solved effectively by the active-set Newton method and the convolutional neural network (CNN), respectively. The proposed ADMM-Newton-CNN approach is validated to be very efficient for the 2-dimensional space case with fine mesh discretization. This work enhances the current literatures in which only the 1-dimensional space case with coarse mesh discretization can be numerically tackled for some smoothing and thus inaccurate surrogate models.

A relevant yet much more challenging problem is to solve TV-regularized models for identifying the diffusion coefficient in an elliptic equation with observation data of function values of the solution, as studied in \cite{ChanT:2003, ChenZ:1999}. This problem is nonconvex and thus intrinsically different from the convex model (\ref{eq:ROTV}). To extend the proposed ADMM-Newton-CNN approach to this nonconvex problem, it is keen to consider how to handle the nonconvex subproblems both theoretically and numerically. It is also interesting to extend the philosophy of algorithmic design, as well as the numerical techniques initiated in this paper, to other parameter identification problems for diffusion coefficients and advection coefficients arising in some elliptic systems with other types of objective functionals, or in some complicated PDE systems.

\begin{table}[htb]
  \centering
  \caption{Parameters $\delta$, $\beta$ and $\sigma$ for Example \ref{exm2}.}\label{parameter_ex2}
    \vspace{3pt}
    \setlength{\tabcolsep}{18pt}\small
    \begin{tabular}{|c|c|c|c|}
	  \hline
      $\delta$ &$\beta$&$\sigma$&$(\beta\sigma)\slash\delta$\\
	  \hline
	  $0.01$ &0.06&12 &84\\
	  \hline
	  $0.05$ &0.3&12&84\\
	  \hline
	    $0.1$ &0.3&24&84\\
	  \hline
  \end{tabular}
\end{table}

\begin{table}[htb]
  \centering
  \caption{Numerical performance for Example \ref{exm2} after the first 50 iterations. }\label{mesh_a_ex2}
  \vspace{3pt}\setlength{\tabcolsep}{5pt}
  \begin{tabular}{|c|c|c|c|c|c|c|}
    \hline
    $\delta$&$h$&Total Newton No.&Total PCG No. for \eqref{eq:A_{k_m}}/\eqref{eq_Hm}&CPU Time (s) & $\|q_h^{50}-q\|_{L^2(\Omega)}/\|q\|_{L^2(\Omega)}$\\
   \hline
    &1/64&   53& 617 / 3156&2.650 & 0.0222\\
       \cline{2-6}
    &1/128& 53&636 /  3543&9.815&0.0136\\
     \cline{2-6}
   0.01 &1/256&53&662 / 3801&46.936& 0.0104\\
 \cline{2-6}
    &1/512& 53&  673 /  3973&193.944& 0.0094\\
     \cline{2-6}
    &1/1024&53&710 / 4059&  841.498& 0.0096\\
    \hline
    \hline
        &1/64&   54&  655 / 1718&2.268& 0.0852\\
       \cline{2-6}
    &1/128&53&656 / 1824&7.964&0.0545\\
     \cline{2-6}
   0.05 &1/256&53&679 / 1932&37.295& 0.0388\\
 \cline{2-6}
    &1/512& 53&  696 / 1985&154.104& 0.0404\\
     \cline{2-6}
    &1/1024&53&722 / 2012&666.666& 0.0604\\
    \hline
    \hline
        &1/64&   54& 651 / 1591&2.278& 0.2160\\
       \cline{2-6}
    &1/128& 54&682 / 1772&8.030&0.1903\\
     \cline{2-6}
   0.1 &1/256&53&690 / 1829&36.437& 0.1558\\
 \cline{2-6}
    &1/512&53&715 / 1915&152.632& 0.1419\\
     \cline{2-6}
    &1/1024&51& 738 / 1967&662.879& 0.1462\\
    \hline
  \end{tabular}
  \end{table}
\begin{table}[htb]
  \centering
  \caption{Computation time for Example \ref{exm2} with $h=1/256$ after the first 50 iterations.}\label{performance_ex2}
  \vspace{3pt}
  \setlength{\tabcolsep}{3pt}\small
  \begin{tabular}{|c|c|c|c|c|c|c|c|c|c|c|c|}
	\cline{1-4} \cline{5-8} \cline{9-12}
	&\multicolumn{3}{c|}{$\delta$=0.01}&&\multicolumn{3}{|c|}{$\delta$=0.05}&&\multicolumn{3}{|c|}{$\delta$=0.1}\\
\cline{1-4} \cline{6-8} \cline{10-12}
	Subtasks&No.&Total time (s)& $\%$Time&&No.&Total time (s)& $\%$Time&&No.&Total time (s)& $\%$Time\\
	\cline{1-4} \cline{6-8} \cline{10-12}
	Newtonian system \eqref{eq_Hm} &53 &19.773&$42.1\%$&       & 53&10.182&$27.3\%$&             &53&9.542&$26.2\%$\\
\cline{1-4} \cline{6-8} \cline{10-12}
    Implementation of CNN &50&9.729&$20.7\%$&                 &50&9.730&$26.1\%$&                     &50&9.425&$25.9\%$\\
\cline{1-4} \cline{6-8} \cline{10-12}
 $N_{k_m}$&53&4.064&$8.7\%$&                &53&4.035 &$10.8\%$&                &53&4.025&$11.0\%$\\
\cline{1-4} \cline{6-8} \cline{10-12}
$A_{k_m}$&53&3.374&$7.2\%$&                &53&3.367&$9.0\%$&                        &53&3.361&$9.2\%$\\
\cline{1-4} \cline{6-8} \cline{10-12}
Linear system \eqref{eq:A_{k_m}} &53&2.792&$5.9\%$&              &53&2.833&$7.6\%$&              &53&2.906&$8.0\%$\\
\cline{1-4} \cline{6-8} \cline{10-12}
Others &&7.205&$15.4\%$&                                              &&7.149&$19.2\%$&                 &&7.178&$19.7\%$\\
\cline{1-4} \cline{6-8} \cline{10-12}
Total&&46.936&$100\%$&                    &&37.295&$100\%$&                    &&36.437&$100\%$\\
	\hline
  \end{tabular}
\end{table}


\phantomsection\small
\addcontentsline{toc}{section}{References}
\bibliographystyle{siam} 


\end{document}